\documentclass[english,a4paper,10pt]{amsart}
\usepackage{babel,amsmath,amssymb,amsbsy,amsfonts,latexsym,amsthm, amscd,amsxtra, amsfonts}
\usepackage{graphicx}
\usepackage[all]{xy}
\usepackage[T1]{fontenc}
\usepackage{graphics,latexsym}

\newtheorem{Teo}{Theorem}[section]
\newtheorem{Def}[Teo]{Definition}
\newtheorem{Cor}[Teo]{Corollary}
\newtheorem{Lem}[Teo]{Lemma}
\newtheorem{Prop}[Teo]{Proposition}

\newtheorem{Rema}[Teo]{Remark}
\newenvironment{Rem}{\begin{Rema} \begin{upshape}} {\end{upshape}\end{Rema}}
\newtheorem*{Pf}{Proof}
\newenvironment{Proof}{\begin{Pf} \begin{upshape}} {\end{upshape} \qed\end{Pf}}

\newcommand\beq[1]{ \begin{equation} \label{#1}}
\newcommand{\eeq}{ \end{equation} }
\newcommand{\beqno}{ \begin{equation*}}
\newcommand{\eeqno}{ \end{equation*}}
\newcommand\beqa[1]{ \begin{eqnarray} \label{#1}}
\newcommand{\eeqa}{ \end{eqnarray} }
\newcommand{\beqano}{ \begin{eqnarray*} }
\newcommand{\eeqano}{ \end{eqnarray*} }

\def\be{\begin{equation}}
\def\ee{\end{equation}}
\def\bea{\begin{eqnarray}}\def\eea{\end{eqnarray}}
\def\nn{\nonumber}

\newcommand{\T}{ {\mathbb T}   }

\newcommand{\R}{ {\mathbb R}   }
\newcommand{\Z}{ {\mathbb Z}   }

\newcommand{\K}{ {\mathbb K}   }

\newcommand{\A}{ {A}   }

\renewcommand \a {\alpha}

\renewcommand \b  {\beta}

\newcommand \m {\mu}

\newcommand \om {\omega}

\newcommand \f {\varphi}

\newcommand \g {\gamma}

\newcommand \s {\sigma}

\renewcommand \L {\Lambda}

\def\cprime{$'$}

\newcommand \cA {{\mathcal A}}
\newcommand \cN {{\mathcal N}}

\newcommand \calM {\mathfrak{M}}

\newcommand \cM {{\mathcal M}}

\newcommand \cL {{\mathcal L}}

\newcommand \dpr {\partial}

\newcommand \rH {{H}}
\newcommand \rT {{T}}

\newcommand \cO {{\mathcal O}}
\newcommand \h {\mathfrak{h}}
\newcommand \Arnold {{Arnol\cprime d }}

\let\0=\noindent

\def\ie{\hbox{\it i.e.,\ }}

\let\dpr=\partial

\let\==\equiv

\newcommand\pb[2] {\{#1, #2\}}

\title[\tiny{A variational approach to the study of the exist. of inv. Lagr. graphs}]{A variational approach to the study of the existence of invariant Lagrangian graphs}
\author{Alfonso Sorrentino}
\address{Dipartimento di Matematica, Universit\`a degli Studi Roma Tre, Largo S. Leonardo Murialdo 1, 00146 Rome (Italy).}
\email{sorrentino@mat.uniroma3.it}
\date{\today}
\subjclass[2000]{37J50, 37J15, 37J35, 37J30}

\begin{document}

\begin{abstract}
This paper surveys some recent results by the author and some collaborators, on the existence of invariant Lagrangian graphs for Tonelli Hamiltonian systems.
\end{abstract}

\maketitle

\tableofcontents


\section{Introduction}\label{sec1}
In this article I would like to describe some properties of Hamiltonian and Lagrangian systems, with particular attention to the relation between their {\it action-minimizing} properties, their {\it symplectic} nature and their {\it dynamics}. More specifically, I shall illustrate what kind of information the {\it principle of least Lagrangian action}\footnote{``{\it Nature is thrifty in all its actions}'', Pierre-Louis Moreau de Maupertuis (1744). A better--known special case of this principle is what is generally called {\it Maupertius' principle}.
Actually, as I learnt from Leo Butler \cite{Butlertalk}, K\"onig published a note claiming priority for Leibniz in the Berlin Academy correspondences overseen by Maupertuis. Priority dispute brought in Euler, Voltaire and ultimately a committee convened by the King of Prussia. In 1913, the Berlin Academy reversed its previous decision and found Leibniz had priority. 
} 
 conveys into the study of  the integrability of these systems, and, more generally, how this information relates to the  existence or to the non--existence of invariant Lagrangian graphs.

Hereafter  I shall address how these very interesting (and difficult) questions can be tackled from different perspectives.
The main results  are all contained in \cite{ButlerSorrentino}, \cite{MassartSorrentino} and \cite{SorrentinoTAMS}, to which I shall refer for a more comprehensive presentation and for detailed proofs.\\

\subsection{Weak Liouville--\Arnold Theorem} \label{sec1.1} (See \cite{SorrentinoTAMS, ButlerSorrentino}). 
It is natural to expect that ``sufficiently'' symmetric systems ought to possess an abundance of invariant Lagrangian graphs. This is indeed the content of a very 
classical result in the study of Hamiltonian systems:  what is generally called {\it Liouville--\Arnold Theorem}  (see for example \cite{Arnoldbook}). This theorem  is concerned with the {\it integrability} of a Hamiltonian system, {\it i.e.}, with the existence of  a regular foliation of the phase space by invariant Lagrangian submanifolds. This theorem provides sufficient conditions for the existence of such a foliation in terms of the existence of  {\it independent} ``symmetries'' that are in {\it involution}.
In order to make the latter condition clearer, let us recall some terminology.
Consider a Hamiltonian system with $n$ degrees of freedom, given by a Hamiltonian $H: V \rightarrow \R$ defined on a $2n$-dimensional symplectic manifold $(V,\omega)$ and denote by $\{\cdot,\cdot\}$ the associated {\it Poisson bracket}, defined as follows: if $f,g \in C^1(M)$, then $\{f,g\}=\om(X_f,X_g)=df\cdot X_g$, where $X_f$ and $X_g$ denote the Hamiltonian vector fields associated to $f$ and $g$ (see for instance \cite{Arnoldbook}).  
An important role in the study of the dynamics  is played by the functions $F: V \rightarrow \R$ that are in  {\it involution} with the Hamiltonian, \ie whose Poisson bracket $\{H,F\}\equiv 0$ on $V$ (equivalently we can say that $H$ and $F$ {\it Poisson-commute}). Such functions, whenever they exist, are called {\it integrals of motion} (or {\it first integrals}) of $H$. It is quite easy to check that the condition of being in involution  is equivalent to asking that $F$ is constant along the orbits of the Hamiltonian flow of $H$ and {\it vice versa}; moreover, this implies that the associated Hamiltonian vector fields $X_H$ and $X_F$ commute.
Liouville--\Arnold theorem relates the integrability of a given Hamiltonian system to the existence of ``enough''  integrals of motion in involution.\\

\noindent{\bf Theorem [Liouville-\Arnold].} {\it
Let $(V,\omega)$ be a symplectic manifold with $\dim V=2n$ and let $H: V \longrightarrow \R$ be a proper Hamiltonian. Suppose that there exists $n$ integrals of motion $F_1, \ldots, F_n : V \longrightarrow   \R$ such that:
\begin{itemize}
\item[{\rm i)}] $F_1,\;\ldots,\; F_n$ are $C^2$ and functionally independent almost everywhere on $V$;
\item[{\rm ii)}] $F_1,\ldots, F_n$ are pairwise in involution, \ie $\pb{F_i}{F_j}=0$ for all $i,j=1,\ldots n$.
\end{itemize}
\noindent Suppose the non-empty regular level set 
$\Lambda_{a} := \{F_1=a_1,\ldots, F_n=a_n\}$ is connected. Then $\Lambda_a$ is an $n$-torus, $\T^n$,
and there is a neighbourhood ${\mathcal O}$ of $0 \in H^1(\Lambda_a;\R)$ such
that for each $c' \in {\mathcal O}$ there is a unique smooth Lagrangian
$\Lambda_{c'}$ that is a graph over $\Lambda_a$ with cohomology class
$c'$. Moreover, the flow of $X_H|\Lambda_{c'}$ is conjugated to a rigid rotation.\\
}

\noindent See for instance \cite[Section 49]{Arnoldbook} for a proof of this theorem.

\begin{Rem}  \label{rem1.1}
The map $F:=(F_1,\ldots, F_n)$ is referred to as an {\it integral
    map}, {\it first-integral map} or a {\it momentum map}.
    The invariance of the level set $\Lambda_a$  simply follows from $F_i$ being integrals of motion; the fact that it is a Lagrangian torus and that the Hamiltonian flow is conjugate to a rigid rotation, strongly relies on these integrals being pairwise in involution and independent. 
 \end{Rem}   
    
A natural question is then the following: {\it Is it possible to weaken the assumptions in the Liouville--\Arnold theorem?} In particular:  {\it What happens when the involution hypothesis on the integrals of motion is dropped?} \\
It is clear from Remark \ref{rem1.1} that the involution hypothesis is fundamental in deducing the non-trivial fact that these invariant sets are Lagrangian and that the motion on them is conjugate to a rigid rotation. Roughly speaking this is not just a sufficient condition, it is also somehow necessary:  the Lagrangianeity of these submanifolds, in fact, is essentially equivalent to the involution hypothesis. 
Hence, it seems almost hopeless to deduce interesting results without assuming it. However, in the case of {\it Tonelli Hamiltonians}   it turns out to be possible.\\
Recall that a Hamiltonian $H\in C^2(T^*M)$ is called {\it Tonelli} if it is fibrewise strictly convex and enjoys fibrewise superlinear growth (see subsection \ref{sec2.2} for a precise definition).
\\

In \cite{SorrentinoTAMS}, I introduced the following definition.

\begin{Def}[{\bf Weak integrability \cite{SorrentinoTAMS}}]
  \label{def:1}
  Let $H \in C^2(T^*M)$. If
  there is a $C^2$ map $F : (T^*M)^n \longrightarrow \R^n$ whose singular set is nowhere
  dense, and $F$ Poisson-commutes with $H$, then we say that $H$ is
{\it weakly integrable}.
\end{Def}

\begin{Rem}
There exist examples of Hamiltonian systems that are weakly integrable but not Liouville Integrable. See for example \cite[Appendix A]{SorrentinoTAMS}. Moreover,
Butler and Paternain \cite{ButlerPaternain} proved the following. Let G be a compact semi-simple Lie group of rank at least 2. In any neighbourhood of the bi-invariant metrics, there are left-invariant metrics with positive topological entropy that are not completely integrable. Nevertheless, these metrics are weakly integrable.
\end{Rem}

In \cite{SorrentinoTAMS} and in a subsequent joint work with Leo Butler \cite{ButlerSorrentino},  we proved a version of Liouville--\Arnold Theorem for weakly integrable Tonelli Hamiltonians. We named this theorem {\it Weak Liouville--\Arnold Theorem}  because it drops the involutivity hypothesis of the classical theorem and still obtains results that are quite analogous. \\

\begin{Teo}[{\bf Weak Liouville-\Arnold}]
\label{mainthm}
Let $M$ be a closed manifold of dimension $n$ and $H:T^*M \longrightarrow \R$ a weakly integrable Tonelli Hamiltonian with integral map
$F: T^*M \longrightarrow \R^n$. If for some cohomology class $c\in H^1(M;\R)$ the Mather set\footnote{A definition of the {\it Mather set} and a description of its properties will be provided in subsection \ref{sec2.4}.}
 $\cM^*_c \subset \; {\rm reg}{F}$, then there exists an open neighbourhood $\cO$ of $c$ in $H^1(M;\R)$ such that the following holds.
\begin{itemize}
\item[{\rm i)}] For each $c'\in \cO$ there exists a smooth invariant Lagrangian graph $\Lambda_{c'}$ of cohomology class $c'$, which admits the structure of a smooth $\T^d$-bundle over a base $B^{n-d}$ that is parallelisable, for some $d>0$. 
\item[{\rm ii)}] The motion on each $\Lambda_{c'}$ is Schwartzman strictly ergodic (see \cite{FGS}), \ie all invariant probability measures have the same rotation vector and the union of their supports equals $\Lambda_{c'}$. In particular, all orbits  are conjugate by a smooth diffeomorphism isotopic to the identity.
\item[{\rm iii)}] Mather's $\alpha$-function (or { minimal average action}) $\alpha: H^1(M;\R) \longrightarrow \R $ is differentiable at all  $c'\in\cO$ and its convex conjugate
$\beta: H_1(M;\R) \longrightarrow \R$ is differentiable at all  rotation vectors $h\in\partial \alpha(\cO)$, where $\partial \alpha(\cO)$ denotes the set of subderivatives of $\alpha$ at some element of $\cO$. In particular, for $c\in \cO$  orbits on $\Lambda_{c'}$ have rotation vector $\partial \alpha(c')$.
\end{itemize}
\end{Teo}

\begin{Rem}
Observe that the flow of $X_H|\Lambda_{c'}$ is a
  	rotation on the $\T^d$--fibres of $\Lambda_{c'}$ with rotation vector
  	$h_{c'}=\partial \alpha(c')$, where $\partial \alpha(c')$ is the
  	derivative of $\alpha$ at $c'$ and $d=b_1(M)$, the first Betti number of $M$. This is analogous to what happens in the classical Liouville-\Arnold
  	theorem, where the rotation vector is the derivative of $H$ (in action--angle coordinates) at $c'$.\\
\end{Rem}

This theorem shows the existence of a family of smooth invariant Lagrangian graphs $\{\Lambda_{c'}\}_{c'\in\cO}$, which form a lamination of the space. However, in general, these graphs cannot be expected to foliate any open set of the phase space. In fact, there are obvious dimensional obstructions: $\{\Lambda_{c'}\}_{c'\in\cO}$ is a $(\dim H^1(M;\R))$--dimensional family of graphs of dimension $\dim M$. Hence, a necessary condition for this to happen is that $\dim H^{1}(M;\R) \geq \dim M$. 
{\it Is this condition the unique obstruction? What can be said if these graphs foliate an open set?}\\
The following result answers these questions (see \cite{SorrentinoTAMS, ButlerSorrentino}).

\begin{Teo}
  \label{maincor}
  Assume the hypotheses of Theorem \ref{mainthm}. 
If $\dim H^1(M;\R) \geq \dim M$, then $\dim H^1(M;\R) = \dim M$
 and $M$ is diffeomorphic to $\T^n=\R^n/\Z^n$. In particular, it follows that $H$ is integrable in the sense of Liouville and therefore the integrals of motion are in involution. \\
 \end{Teo}
  
  \begin{Rem} ({\it i}) In other words, if $\dim H^1(M;\R) \geq \dim M$ then the notion of { weak integrability} is equivalent to the classical notion of  {Liouville integrability}. This hypothesis is satisfied, for example, if $M=\T^n$.
  
({\it ii})  Actually, in \cite{ButlerSorrentino} we proved something more (these results go beyond the purpose of this survey article, therefore we refer all interested readers  to \cite{ButlerSorrentino} for complete proofs). 
 We showed that  under the hypotheses of Theorem \ref{mainthm}:  
\begin{itemize}
\item[a)] if $\dim M \leq 3$, then $M$ is diffeomorphic to a torus; 
\item[b)] if $\dim M=4$, then -- assuming the virtual Haken conjecture\footnote{In $3$-manifold topology, a central role is played by those closed
$3$-manifolds which contain a non-separating {\it incompressible}
surface, or dually, which have non-vanishing first Betti number. Such
manifolds are called {\it Haken}; it is an outstanding conjecture
that every irreducible $3$-manifold with infinite fundamental group
has a finite covering that is Haken \cite[Questions
1.1--1.3]{MR739142}. This conjecture is implied by the 
virtually fibred conjecture \cite{MR2399130}. Given the proof of the
geometrisation conjecture, the virtual Haken conjecture is proven for
all cases but hyperbolic $3$-manifolds. Thurston and Dunfield have
shown there is good reason to believe the conjecture is true in this
case \cite{MR1988291}.
}
--   $M$ is diffeomorphic to either $\T^4$ or $\T^1 \times E$,  where $E$ is an orientable $3$-manifold finitely covered by      $S^3$.
\end{itemize}
See \cite[Theorem 1.2]{ButlerSorrentino} for a detailed proof.

({\it iii}) Moreover, we investigated the case in which the system's
symmetries are not classical and do not come from conserved
quantities, but are induced by invariance under the action of an
amenable Lie group on the universal cover of the manifold. This action
need not descend to the quotient and is generally only evident in
statistical properties of orbits. In particular, these symmetries may
only manifest themselves in the structure of the action-minimizing
sets. Analogous results to the ones discussed above can be also proven in this case; see \cite[Theorem 1.3]{ButlerSorrentino} for a precise statement.\\
\end{Rem}


\subsection{Differentiability of  the minimal average action and Integrability} \label{sec1.2} (See \cite{MassartSorrentino}).
In this section I would like to discuss another possible approach to the study of the existence of invariant Lagrangian graphs and the integrability of the system.

  In the study of Tonelli Lagrangian and Hamiltonian systems,  a central role in understanding  the dynamical and topological properties of the  action-minimizing sets (see subsection \ref{sec2.4}), is played by the so-called {\it Mather's average action} (sometimes referred to as $\beta$-{\it function}  or {\it effective Lagrangian}), with particular attention to its differentiability and non-differentiability properties. 
Roughly speaking,  this is a convex superlinear function on the first homology group of the base manifold, which represents the minimal action of invariant probability measures within a prescribed  {\it homology class}, or {\it rotation vector} (see (\ref{defbeta}) for a more precise definition). 
Understanding whether or not this function is differentiable, or even smoother, 
and what are the implications of its regularity to the dynamics  of the system is a formidable problem, which is still far from being understood.
Examples of Lagrangians  admitting a smooth $\beta$-function are easy to construct; trivially, if the base manifold $M$ is such that $\dim \rH_1(M;\R)=0$ then $\beta$ is a function defined on a single-point set and it is therefore smooth. Furthermore, if $\dim \rH_1(M;\R)=1$ then a result by Carneiro \cite{Carneiro} (see also Lemma \ref{Carneiro radial}) allows one to conclude that $\beta$ is differentiable everywhere, except possibly at the origin. 
As soon as  $\dim \rH_1(M;\R)\geq 2$ the situation becomes definitely less clear and the smoothness of $\beta$ becomes a more ``untypical'' phenomenon. Nevertheless, it is still possible to find some interesting examples in which it is smooth.
 For instance, let 
$H:\rT^*\T^n \longrightarrow \R$ be a completely integrable Tonelli Hamiltonian system, given by $H(x,p)=\h(p)$, and consider the associated Lagrangian $L(x,v)=\ell(v)$ on $\rT \T^n$.  It is easy to check (see subsection \ref{sec2.3}) that in this case, up to identifying $\rH_1(\T^n;\R)$ with $\R^n$, one has $\beta(h)=\ell(h)$ and therefore $\b$ is as smooth as the Lagrangian.  One can weaken the assumption on the complete integrability of the system and consider $C^0$-{\it integrable systems}, \ie Hamiltonian systems that admit a foliation of the phase space by disjoint invariant continuous Lagrangian graphs, one for each possible cohomology class (see Definition \ref{C0integrability} and \cite{Arnaud}). It is then possible to prove that also in this case the associated $\beta$ function is $C^1$.

These observations raise the following question: {\it If $\dim H_1(M;\R)\geq 2$, does the regularity of $\b$ imply the integrability of the system?} Or more generally, {\it is the existence of an invariant Lagrangian graph detected by some regularity property of this function?}\\

In a joint work with Daniel Massart \cite{MassartSorrentino}, we addressed the above problem in the case of Tonelli Lagrangians on closed surfaces, not necessarily orientable
(in this latter case, one considers the lifted Lagrangian to the orientable double cover).
The main results can be summarized as follows.

\begin{Teo}\label{teoms}
Let $M$ be a closed surface and $L:\rT M \longrightarrow \R$ a Tonelli Lagrangian. 
\begin{itemize}
\item[{\rm (i)}]   If $M$ is  not the sphere,  the projective plane, the Klein bottle or the torus, then $\beta$ cannot be $C^1$ on the whole of $\rH_1(M;\R)$.
\item[{\rm (ii)}] If $M$ is  not the torus then the Lagrangian cannot be $C^0$-integrable.
\item[{\rm (iii)}]	If $M$ is the torus, then $\beta$ is $C^1$ if and only if the system is $C^0$-integrable.
\end{itemize}
\end{Teo}

In particular, in the orientable case it is possible to relate the differentiability of $\beta$ at (non-singular\footnote{We refer to Section \ref{sec4} and \cite{MassartSorrentino} for a more precise definition}) $1$-irrational homology classes to the existence of invariant Lagrangian graphs foliated by periodic orbits. 

\begin{Teo}
Let  $M$ be a closed, oriented  surface and  $L: T M \longrightarrow \R$ an autonomous  Tonelli Lagrangian. Let $h_0$ be a 1-irrational, non-singular homology class, and let
$\langle\partial \beta (h_0)\rangle $ denote 
the underlying vector space of the affine subspace generated by $\partial \beta(h_0)$  in $H^1(M,R)$.
Then, the dimension of $\langle \partial \beta (h_0) \rangle$ is at least the genus of $M$. Moreover, if $M=\T^2$ and $\beta$ is differentiable at $h_0$, then there exists an invariant Lagrangian graph foliated by periodic orbits, whose homology is a multiple of $h_0$.
\end{Teo}

\begin{Rem}
Observe that the above result is not true if $h$ is singular (see example \cite[Remark 2]{MassartSorrentino}). Moreover, 
when $M$ is not orientable, the situation is different because $\beta$ may have flats of maximal dimension (that is, of dimension equal to the first Betti number of $M$). So $\beta$ may well be differentiable at some 1-irrational, non-singular homology class without having any invariant Lagrangian graph foliated by periodic orbits (see subsection \ref{sec4.2} and \cite[pages 11-12]{MassartSorrentino} for a more precise discussion of this issue).
\end{Rem}

The above results consider the case of $C^0$--integrable case. An open question is whether similar results can be obtained with $C^0$--integrability replaced by integrability in the sense of Liouville.
In the case of mechanical systems  we can bridge this gap. A mechanical Lagrangian is a Lagrangian of the form
$L(x,v) = 1/2 \,  g_x(v,v) +f(x)$, where $g$ is a Riemannian metric and $f$ is a $C^2$ function on $\T^2$  (see also subsection \ref{sec2.2} for the definition).

\begin{Prop}\label{propmech}
Let $L$ be a  mechanical Lagrangian on a $2$-dimensional torus, whose $\beta$-function is $C^1$. Then, the potential $f$ is identically constant and the metric $g$ is flat. In particular, $L$ is integrable in the sense of Liouville.
\end{Prop}

\noindent See \cite[Proposition 6]{MassartSorrentino} for a complete proof.

\subsection{Outline of the article}
In order to make the material interesting and accessible to a wider audience, in Section \ref{sec2} I shall provide a brief introduction to Mather's theory, starting from an historical 
excursion (subsection \ref{sec2.1}) and guiding the  reader through a cartoon example (subsection \ref{sec2.2}). A description of the general theory will be provided in subsection \ref{sec2.3}.
In section \ref{sec3}, I shall describe the results related to weakly integrable system, providing a sketch of the main ideas involved in the proofs. Finally, in section \ref{sec4}
I shall survey the results about the differentiability of the minimal average action and its relation the the integrability of the system.
Due to the nature of this survey, most of the proofs will not be provided. However, I shall refer the interested reader to the relevant literature.\\

\noindent{\Small {\rm ACKNOWLEDGEMENTS}}. The author is grateful to {\it Unione Matematica Italiana} (UMI) for the invitation to deliver this talk at 
{\it XIX Congresso dell'UMI} (session on {\it Nonlinear Analysis and dynamical systems}), held in Bologna (Italy) from 12th to 17th September 2011. \\


\section{Introduction to (Aubry-)Mather theory}\label{sec2}

\subsection{From KAM theory to AM theory}\label{sec2.1}

The celebrated Kolmogorov-\Arnold-Moser (or KAM) theorem finally settled the old question concerning the existence of {\it quasi-periodic}  motions for {\it nearly-integrable} Hamiltonian systems, \ie Hamiltonian systems that are slight perturbation of an {integrable} one.  In the integrable case, in fact, the whole phase space is foliated by invariant {Lagrangian} submanifolds that are diffeomorphic to tori, and on which the dynamics is conjugate to a rigid rotation. On the other hand, it is natural to ask  what happens to such a foliation and to these {\it stable motions} once the system is perturbed.
In 1954 Kolmogorov \cite{KolmogorovKam} - and later \Arnold  \cite{ArnoldKAM} and Moser \cite{MoserKAM} in different contexts - proved that, in spite of the generic disappearence of the invariant 
tori filled by periodic orbits, already pointed out by Henri Poincar\'e,  for small perturbations of an integrable system it is still possible to find invariant Lagrangian tori corresponding to  ``{\it strongly non-resonant}''  rotation vectors. This result, commonly referred to as {\it KAM theorem}, from the initials of the three main pioneers, besides opening the way to a new understanding of the nature of Hamiltonian systems and their stable motions, contributed to raise new interesting questions, for instance: {\it what about the destiny of the {\it stable motions} (orbits on KAM tori) that are destroyed by effect of the perturbation? Is it possible to identify something reminiscent of their past presence? What can be said about a system which is not close to an integrable one?}

{\it Aubry-Mather theory} provides an answer to these questions.  Developed independently by Serge Aubry \cite{Aubry} and 
John Mather \cite{Math82} in the eighties, this novel approach to the study of the dynamics of
{\it twist diffeomorphisms of the annulus} (which correspond to Poincar\'e maps of $1$-dimensional non-autonomous Hamiltonian systems) pointed out the existence of many {\it action-minimizing sets}, which in some sense generalize invariant rotational curves and that always exist, even after rotational curves are destroyed. Besides providing a detailed structure theory for these new sets, this powerful approach yielded  a better understanding of the destiny of invariant rotational curves and to the construction of interesting chaotic orbits as a result of their destruction \cite{Mather86, Matherorbitsdiffeo}.

Motivated by these achievements, John Mather \cite{Mather91, Mather93} - and later Ricardo Ma\~n\'e \cite{ManeI, Maneminmeasure}
and Albert Fathi \cite{Fathibook} in different ways - developed a generalization of this theory to higher dimensional systems.
Positive definite superlinear Lagrangians on compact manifolds, also called {\it Tonelli Lagrangians} (see Definition \ref{defTonelliLag}),
were the appropriate setting to work in.  
Under these conditions, in fact, it is possible to prove the existence of interesting invariant  sets, known as {\it Mather}, {\it Aubry} and {\it Ma\~n\'e} sets, which generalize KAM tori, and which continue to exist even after KAM tori's disappearance or when it does not make sense to speak of them (for example when the system is ``far'' from any integrable one). 

In the following we shall provide a brief introduction to Mather's theory. We shall first discuss a cartoon example (Integrable systems) and then show how similar ideas can be extended to a more general setting. A comprehensive exposition of this material (and much more) can be also found in \cite{SorrentinoLectNotes} (see also \cite{bangert, Matherforni, Mather10, Fathibook} for other interesting related material).


\subsection{Tonelli Lagrangians and Hamiltonians on compact manifolds}\label{sec2.2}
Before starting, let us introduce the basic setting that we shall consider in the following.
Let $M$ be a compact and connected smooth manifold without boundary.
Denote by $TM$ its tangent bundle and $T^*M$ the cotangent one. A
point of $TM$ will be denoted by $(x,v)$, where $x\in M$ and $v\in
T_xM$, and a point of $T^*M$ by $(x,p)$, where $p\in T_x^*M$ is a
linear form on the vector space $T_xM$. Let us fix a Riemannian
metric $g$ on it and denote by $d$ the induced metric on $M$; let
$\|\cdot\|_x$ be the norm induced by $g$ on $T_xM$; we shall use the same notation for the norm induced on
$T_x^*M$.

We shall consider functions $L:\rT M \longrightarrow \R$ of class $C^2$, which are called {\it Lagrangians}. Associated to each Lagrangian, there is a flow on $\rT M$ called the {\it Euler-Lagrange flow}, defined as follows. Let us consider the action functional $\A_L$ from the space of continuous
piecewise $C^1$ curves $\g:[a,b]\rightarrow M$, with $a\leq b$, defined by:
$$
\A_L(\g) := \int_a^b L(\g(t),\dot{\g}(t))\,dt.
$$
Curves that extremize\footnote{These extremals are not in general minima. The existence of global minima and the study of the corresponding motions is the core of Aubry-Mather theory; see subsection \ref{sec2.4}} this functional among all curves with the same end-points (and the same time-length) are solutions of the {\it Euler-Lagrange equation}:
\beqa{ELequation}
\frac{d}{dt}\frac{\dpr L}{\dpr v} (\g(t),\dot{\g}(t)) = \frac{\dpr L}{\dpr x} (\g(t),\dot{\g}(t)) \qquad \forall\,t\in[a,b]\,.
\eeqa
Observe that this equation is equivalent to
$$
\frac{\dpr^2 L}{\dpr v^2} (\g(t),\dot{\g}(t)) \ddot{\g}(t) = \frac{\dpr L}{\dpr x} (\g(t),\dot{\g}(t)) - \frac{\dpr^2 L}{\dpr v \dpr x} (\g(t),\dot{\g}(t))\dot{\g}(t)\,,$$
therefore, if the second partial vertical derivative ${\dpr^2 L}/{\dpr v^2}(x,v)$ is non-degenerate at all points of $\rT M$, we can solve for $\ddot{\g}(t)$. This condition
$$
\det \frac{\dpr^2 L}{\dpr v^2} \neq 0
$$
is called {\it Legendre condition} and allows one to define a vector field $X_L$ on $\rT M$, such that the solutions of 
$\ddot{\g}(t)=X_L(\g(t),\dot{\g}(t))$ are precisely the curves satisfying the Euler-Lagrange equation.
This vector field $X_L$ is called the {\it Euler-Lagrange vector field} and its flow ${\Phi^L_t}$ is the {\it Euler-Lagrange flow} associated to $L$. It turns out that $\Phi^L_t$ is $C^1$ even if $L$ is only $C^2$ (see Remark \ref{ELC1}).

\begin{Def}[{\bf Tonelli Lagrangian}]\label{defTonelliLag}
A function $L:\,TM\, \longrightarrow \,\R$ is called a {\it Tonelli
Lagrangian} if:
\begin{itemize}
\item[i)]   $L\in C^2(TM)$;
\item[ii)]  $L$ is strictly convex in the fibres, in the $C^2$ sense, \ie the second partial vertical derivative
${\dpr^2 L}/{\dpr v^2}(x,v)$ is positive definite, as a quadratic form, for all $(x,v)$;
\item[iii)] $L$ is superlinear in each fibre, \ie
            $$\lim_{\|v\|_x\rightarrow +\infty} \frac{L(x,v)}{\|v\|_x} = + \infty.$$
            This condition is equivalent to ask that for each $A\in \R$ there exists $B(A)\in\R$ such that
            $$ L(x,v) \geq A\|v\| - B(A) \qquad \forall\,(x,v)\in \rT M\,.$$ 
\end{itemize}
\end{Def}
Observe that since the manifold is compact, then condition {\it iii)} is independent of the choice of the Riemannian metric $g$.\\

\noindent {\bf Examples of Tonelli Lagrangians.}
\begin{itemize}
\item {\bf Riemannian Lagrangians.} Given a Riemannian metric $g$ on $\rT M$, the {\it Riemannian Lagrangian} on $(M,g)$ is given by the  {\it kinetic energy}:
$$
L(x,v) = \frac{1}{2} \| v\|_x^2\,.
$$
Its Euler-Lagrange equation is the equation of the geodesics of $g$: 
$$\frac{D}{dt}\dot{x} \equiv 0\,,$$
and its Euler-Lagrange flow coincides with the geodesic flow.
\item {\bf Mechanical Lagrangians.} These Lagrangians play a key-role in the study of classical mechanics. They are given by the sum of the kinetic energy and a {\it potential} $U: M \longrightarrow \R$:
$$
L(x,v) = \frac{1}{2} \| v\|_x^2 + U(x)\,.
$$
The associated Euler-Lagrange equation is given by:
$$\frac{D}{dt}\dot{x} = \nabla U(x)\,,$$
where $\nabla U$ is the gradient of $U$ with respect to the Riemannian metric $g$, \ie
$$
d_x U \cdot v = \langle \nabla U(x), v \rangle_x \quad \forall\, (x,v)\in \rT M\,.
$$
\item {\bf Ma\~n\'e's Lagrangians.} This is a particular class of Tonelli Lagrangians, introduced by Ricardo Ma\~n\'e in \cite{Maneminmeasure}. If $X$ is a $C^k$ vector field on $M$, with $k\geq 2$, one can embed its flow 
$\varphi^X_t$ into the Euler-Lagrange flow associated to a certain Lagrangian, namely 
$$ L_X(x,v)= \frac{1}{2}\left\| v- X(x) \right\|_x^2\,.$$ 
It is quite easy to check that the integral curves of the vector field $X$ are solutions of the Euler-Lagrange equation. In particular, the Euler-Lagrange flow $\Phi^{L_X}_t$ restricted to ${\rm Graph}(X)=\{(x,X(x)),\;x\in M\}$ (that is clearly invariant) is conjugate to the flow of $X$ on $M$
and the conjugacy is given by $\pi|{{\rm Graph}(X)}$,
where $\pi: \rT M \rightarrow M$ is the canonical projection. In other words, the following diagram commutes:
$$\xymatrix{
{{\rm Graph}(X)} \ar@{->}[d]_{\pi} \ar@{->}[r]^{\Phi^{L_X}_t}  & {{\rm Graph}(X)}\ar@{->}[d]^{\pi} 
\\  
{M} \ar@{->}[r]_{\varphi^X_t}  & {M} \\}
 $$
 that is, for every $x\in M$ and every $t\in \R$, $\Phi^{L_X}_t(x,X(x))=(\g_x^X(t),\dot{\g}^X_x(t))$, where $\g_x^X(t)=\varphi_t^X(x)$.\\ 
\end{itemize}

In the study of classical dynamics, it turns often very useful to consider the associated {\it Hamiltonian system}, which is defined  on the cotangent bundle $\rT^* M$. 
Given a Lagrangian $L$, we can define the associated {\it Hamiltonian},
as its Fenchel transform (or {\it Fenchel-Legendre transform}), see \cite{Rockafellar}: 
\beqano H:\; T^*M &\longrightarrow & \R \\
(x,p) &\longmapsto & \sup_{v\in T_xM} \{\langle p,\,v \rangle_x -
L(x,v)\}\, \eeqano where $\langle \,\cdot,\,\cdot\, \rangle_x$
denotes the canonical pairing between the tangent and cotangent
bundles.

If $L$ is a Tonelli Lagrangian, one can easily prove that $H$ is finite everywhere (as a consequence of the superlinearity of $L$), superlinear and strictly convex in each fibre (in the $C^2$ sense). 
Observe that $H$ is also $C^2$. In fact
the Euler-Lagrange vector field corresponds, under the Legendre transformation, to a vector field on $T^*M$ given by Hamilton's equation; it is easily seen that this vector field is $C^1$ (see \cite[p. 207]{caratheodory}).
Such  a Hamiltonian is called {\it Tonelli} (or {\it optical}) {\it Hamiltonian}.

\begin{Def}[{\bf Tonelli Hamiltonian}]\label{opticalHamiltonian}
A function $H:\,T^*M\longrightarrow \R$ is called a 
{\it Tonelli (or optical) Hamiltonian} if:
\begin{itemize}
\item[i)]  $H$ is of class $C^2$;
\item[ii)]  $H$ is strictly convex in each fibre in the $C^2$ 
sense, \ie the second partial 
vertical derivative ${\dpr^2 H}/{\dpr p^2}(x,p)$ is positive definite, 
as a quadratic form, for any $(x,p)\in T^*M$;
\item[iii)] $H$ is superlinear in each fibre, \ie
$$\lim_{\|p\|_x\rightarrow +\infty} \frac{H(x,p)}{\|p\|_x} = + \infty\,.$$
\end{itemize}
\end{Def}

\vspace{20 pt}

\noindent {\bf Examples of Tonelli Hamiltonians.}\\
Let us see what are the Hamiltonians associated to the Tonelli Lagrangians that we have introduced in the previous examples.

\begin{itemize}
\item {\bf Riemannian Hamiltonians.} If $L(x,v) = \frac{1}{2} \| v\|_x^2$ is the Riemannian Lagrangian associated to a Riemannian metric $g$ on $M$, the corresponding Hamiltonian will be
$$
H(x,p) = \frac{1}{2} \| p\|_x^2,
$$
where $\|\cdot\|$ represents - in this last expression - the induced norm on the cotangent bundle $\rT^* M$. 
\item {\bf Mechanical Hamiltonians.} If $L(x,v) = \frac{1}{2} \| v\|_x^2 + U(x)$ is a mechanical Lagrangian, the associated Hamiltonian is:
$$
H(x,p) = \frac{1}{2} \| p \|_x^2 - U(x).
$$
It is sometimes referred to as {\it mechanical energy}.
\item {\bf Ma\~n\'e's Hamiltonians.}  If $X$ is a $C^k$ vector field on $M$, with $k\geq 2$, and 
$ L_X(x,v)= \left\| v- X(x) \right\|_x^2$ is the associated Ma\~n\'e Lagrangian, one can check that the corresponding Hamiltonian is given by:
$$
H(x,p)= \frac{1}{2}\|p\|_x^2 + \langle p,X(x) \rangle\,.
$$
\end{itemize}

\vspace{10 pt}

Given a Hamiltonian one can consider the associated {\it Hamiltonian flow} $\Phi^H_t$ on $\rT^* M$. In local coordinates, this flow can be expressed in terms of the so-called {\it Hamilton's equations}:
\beqa{Hamiltonequation}
\left\{
\begin{array}{l}
\dot{x}(t) = \frac{\dpr H}{\dpr p} (x(t),p(t))\\
\dot{p}(t) = - \frac{\dpr H}{\dpr x} (x(t),p(t))\,.
\end{array}\right.
\eeqa

We shall denote by $X_H(x,p):=\left(\frac{\dpr H}{\dpr p} (x,p),- \frac{\dpr H}{\dpr x} (x,p)\right)$ the {\it Hamiltonian vector field} associated to $H$. This has a more intrinsic (geometric) definition in terms of the canonical symplectic structure $\omega$ on $\rT^*M$. 
In fact, $X_H$ is the unique vector field that satisfies
$$
\omega \left( X_H (x,p), \cdot \right) = d_x H(\cdot)\qquad \forall (x,p)\in\rT^*M.
$$
For this reason, it is sometime called {\it symplectic gradient of} $H$.
It is easy to check from both definitions that - only in the autonomous case - the Hamiltonian is a {\it prime integral of the motion}, \ie it is constant along the solutions of these equations. \\

Now, we would like to explain what is the relation between the Euler-Lagrange flow and the Hamiltonian one. 
It follows easily from the definition of Hamiltonian (and Fenchel transform) that for each $(x,v) \in \rT M$ and $(x,p)\in \rT^* M$ the following inequality holds:
\beqa{Fenchelineq}
 \langle p,\,v \rangle_x \leq L(x,v) + H(x,p)\,.
\eeqa
This is called {\it Fenchel inequality} (or {\it Fenchel-Legendre inequality}, see \cite{Rockafellar}) and it plays a crucial role in the study of Lagrangian and Hamiltonian dynamics and in the variational methods that we are going to describe. In particular, equality holds if and only if $p = \dpr L/\dpr v (x,v)$. One can therefore introduce the following diffeomorphism between $TM$ and $T^*M$, known as  {\it Legendre transform}:
\beqa{Legendretransform} \cL:\; TM &\longrightarrow & T^*M\nonumber\\
(x,v) &\longmapsto & \left(x,\,\frac{\dpr L }{\dpr v}(x,v)
\right)\,. \eeqa
Moreover, the following relation with the Hamiltonian holds:
$$ H \circ \cL(x,v) = \left\langle \frac{\dpr L }{\dpr v}(x,v),\,v \right\rangle_x - L(x,v)\,.$$
A crucial observation is that this diffeomorphism $\cL$ represents a conjugacy between
the two flows, namely the Euler-Lagrange flow on $\rT M$ and the Hamiltonian flow on $\rT^* M$; in other words, the following diagram commutes:
$$\xymatrix{
{\rT M} \ar@{->}[d]_{\cL} \ar@{->}[r]^{\Phi^L_t}  & {\rT M}\ar@{->}[d]^{\cL} 
\\  
{\rT^* M} \ar@{->}[r]_{\Phi^H_t}  & {\rT^* M} \\}
 $$
 
 \begin{Rem}\label{ELC1}
 Since $\cL$ and the Hamiltonian flow $\Phi^H$ are both $C^1$, then it follows from the commutative diagram above that the Euler-Lagrange flow is also $C^1$.\\
\end{Rem}

\vspace{20 pt}



\subsection{Cartoon example: action-minimizing properties of integrable systems}\label{sec2.3}

Before entering into the details of Mather's work, I would like to discuss  a very easy case: properties of invariant measures of an integrable system.
 This will provide us with a better understanding of the ideas behind Mather's theory and will make it clearer in which sense these {\it action-minimizing sets} -- namely, what we shall call {\it Mather sets} (see subsection \ref{sec2.4}) -- represent a generalization of KAM tori.

Let $H: \rT^*\T^n \longrightarrow \R$ be an integrable Tonelli Hamiltonian in action-angle coordinates, {\it i.e.}, $H(x,p)=\h( p)$ depends only on the action variables\footnote{In general one can find these coordinates only locally. For the sake of simplicity,  let us assume in this example that they are defined globally. This simplification will not affect the main purpose of this section}.
On the other hand the associated Tonelli Lagrangian  $L: \rT\T^n \longrightarrow \R$ has the form $L(x,v)=\ell(v)$.
Let us denote by $\Phi^\h$ and $\Phi^{\ell}$ the respective flows and identify $\rT^*\T^n$ and $\rT \T^n$ with $\T^n\times \R^n$. \\

The Hamiltonian flow in this case is very easy to study. Hamilton's equations are:
$$
\left\{
\begin{array}{l}
\dot{x} = \frac{\partial \h}{\partial p}( p) = :\rho(p)\\
\dot{p} = - \frac{\partial \h}{\partial x}( p) = 0,
\end{array}
\right.
$$
therefore $\Phi^\h_t(x_0,p_0) = \left(x_0+t  \rho( p_0),\,p_0 \right)$. In particular, $p$ is an integral of motion, that is, it remains constant along the orbits. Therefore, the phase space $\rT^*\T^n$ is foliated by invariant tori $\Lambda^*_{p_0} = \T^n \times \{p_0\}$ on which the motion is a rigid rotation with rotation vector $\rho( p_0)$. Clearly these tori are Lagrangian.\\

We would like to understand if these invariant tori can be characterized in a different way, possibly in a form that can be more feasibly extended to the non-integrable case.

\begin{Rem}
Observe that on the Lagrangian side, the counterpart of these invariant tori $\Lambda^*_{p_0} \subset \rT^*\T^n$ is given by 
$\tilde{\Lambda}_{p_0} = \T^n \times \{\rho(p_0)\} \subset \rT \T^n$.  It is trivial to deduce that they are also invariant and that they also foliate the whole tangent bundle.
\end{Rem}

We have recalled in Section \ref{sec2.2} that the Euler-Lagrange flow can be equivalently defined in terms of a variational principle associated to the {\it Lagrangian action functional} $A_L$. We would like to study action-minimizing properties of these invariant manifolds; for this, it is much better to work in the Lagrangian setting. Moreover, instead of considering properties of single orbits, it would be more convenient to study ``collection'' of orbits, in the form of {\it invariant probability measures}\footnote{Actually, it is also possible study directly orbits. See Remark \ref{rem2.17}} and consider their {\it average action}. If $\mu$ is an invariant probability measure for $\Phi^L$ -- {\it i.e.}, ${(\Phi^{\ell}_t)^*\mu=\mu}$ for all $t\in \R$, where $(\Phi^{\ell}_t)^*\mu$ denotes the pull-back of the measure -- then we define:
$$
A_{\ell}(\mu):= \int_{\rT\T^n} {\ell}(v)\,d\mu.
$$

Let us consider  any invariant probability measure $\mu_{0}$ supported on $\tilde{\Lambda}_{p_0}$ and compute its action. Observe that on the support of this measure $\ell(v)\equiv \ell(\rho(p_0))$. Then:
\begin{eqnarray}\label{Lactionmu0}
A_{\ell}(\mu_0) &=& \int_{\rT\T^n} \ell(v)\,d\mu_0 = \int_{\rT\T^n} \ell(\rho(p_0))\,d\mu_0 =\nonumber\\
&=& \ell(p_0) =  p_0 \cdot \rho(p_0) - \h(p_0),
\end{eqnarray}
where in the last step we have used the Fenchel-Legendre duality between the $h$ and $\ell$.

Let us now consider a general invariant probability measure $\mu$. In this case it is not true anymore that $\ell(v)$ is constant on the support of $\mu$. However, using Fenchel-Legendre inequality (see (\ref{Fenchelineq})), we can conclude that $\ell(v) \geq  p_0 \cdot v - \h(p_0)$ for each $v\in \R^n$. Hence:

\begin{eqnarray}\label{Lactionmu}
A_{\ell}(\mu) &=& \int_{\rT\T^n} \ell(v)\,d\mu \geq \int_{\rT\T^n}  \left( p_0 \cdot v - \h(p_0) \right) \,d\mu =\nonumber\\
&=& \int_{\rT\T^n}   p_0\cdot v \, d\mu - \h(p_0) =  p_0\cdot \left(  \int_{\rT\T^n}   v \, d\mu  \right) - \h(p_0).
\end{eqnarray}

We would like to compare  expressions (\ref{Lactionmu0}) and (\ref{Lactionmu}). However, in the case of a general measure, we do not know how to evaluate the term 
$\int_{\rT\T^n}   v \, d\mu$. One possible trick to overcome this problem is the following: instead of considering the action of $\ell(v)$, let us consider the action of
$\ell(v) - p_0\cdot v$. It is easy to see that this ``new'' Lagrangian is also Tonelli (we have subtracted a linear term in $v$) and that it has the same Euler-Lagrange flow as $\ell$! In this way we obtain from (\ref{Lactionmu0}) and (\ref{Lactionmu}) that:

\begin{eqnarray*}
A_{\ell - p_0\cdot v}(\mu_0) &=&  - \h(p_0) \qquad {\rm and }\qquad
A_{\ell - p_0\cdot v}(\mu) \geq - \h(p_0),
\end{eqnarray*}

which are now comparable. Hence, we have just showed the following fact:\\

\noindent {\bf Fact 1:} {\it Every invariant probability measure supported on $\tilde{\Lambda}_{p_0}$ minimizes the action $A_{\ell-p_0\cdot v}$ amongst all invariant probability measures of $\Phi^{\ell}$.\\
}

In particular, we can characterize our invariant tori in a different way:
\begin{equation}
\tilde{\Lambda}_{p_0} = \bigcup\{{\rm supp}\, \mu: \; \mu\; \mbox{minimizes} \;  A_{\ell-p_0\cdot v}\}.
\end{equation}

Moreover, there is a relation between the energy (Hamiltonian) of the invariant torus and the minimal action of its invariant probability measures:
\begin{equation}\label{minactionhint}
\min \{A_{\ell-p_0\cdot v}(\mu):\; \mu\; \mbox{is an inv. prob. measure}\} = -\h(p_0).
\end{equation}

\vspace{20 pt}

Observe that it was somehow expectable that we needed to modify the Lagrangian in order to obtain information on a specific invariant torus. In fact, in the case of an integrable system we have a foliation of the space made by these invariant tori and it  would be hopeless to expect that they could all be obtained as extremals of the same action functional.  In other words, what we did was to add a {\it weighting term} to our Lagrangian, in order to magnify some motions rather than others.\\
Is it possible to distinguish these motions in a different way? Let us go back to (\ref{Lactionmu0}) and (\ref{Lactionmu}). The main problem in comparing these two expression was represented by the term $\int_{\rT\T^n}   v \, d\mu$. This can be interpreted as a sort of average rotation vector of orbits in the support of $\mu$. Hence, let us define the {\it average rotation vector of } $\mu$ as:
$$
\rho(\mu):=\int_{\rT\T^n}   v \, d\mu \in \R^n.
$$
We shall give a more precise definition of it (which is also meaningful on manifolds different from the torus) in subsection \ref{sec2.4}. 

Let now $\mu$ be an invariant probability measure of $\Phi^{\ell}$ with rotation vector $\rho(\mu)=\rho(p_0)$. It follows from (\ref{Lactionmu}) that:
\begin{eqnarray*}
A_{\ell}(\mu) &\geq&  p_0\cdot \left(  \int_{\rT\T^n}   v \, d\mu  \right) - \h(p_0) =  p_0\cdot \rho(\mu) - \h(p_0) =\\
&=&p_0\cdot \rho(p_0) - \h(p_0) = \ell(\rho(p_0)).
\end{eqnarray*}

Therefore, comparing with (\ref{Lactionmu0}) we obtain another characterization of $\mu_0$:\\

\noindent {\bf Fact 2:} {\it Every invariant probability measure supported on $\tilde{\Lambda}_{p_0}$ minimizes the action $A_{\ell}$ amongst all invariant probability measures of $\Phi^{\ell}$ with rotation vector $\rho(p_0)$.\\
}

In particular:
\begin{equation}
\tilde{\Lambda}_{p_0} = \bigcup\{{\rm supp}\, \mu: \; \mu\; \mbox{minimizes} \;  A_{\ell} \; \mbox{amongst measures with rot. vect.}\; \rho(p_0)\}.
\end{equation}

Moreover, there is a relation between the value of the Lagrangian at $\rho(p_0)$ and the minimal action of all invariant probability measures with rotation vector $\rho(p_0)$:
\begin{equation}\label{minactionh}
\min \{A_{\ell}(\mu):\; \mu\; \mbox{is an inv. prob. meas. with rot. vect.}\; \rho(p_0)\} = \ell(\rho(p_0)).
\end{equation}

\vspace{10 pt}

\begin{Rem}
One could also study directly orbits on these tori and try to show that their action minimizes a modified Lagrangian action, in the same spirit as we already saw for measures. See \cite{SorrentinoLectNotes} and Remark \ref{rem2.17} for more details.
\end{Rem}

\vspace{20 pt}


\subsection{Mather's theory for Tonelli Lagrangian systems}\label{sec2.4}

In this section  we would like to discuss Mather's theory  for general Tonelli Lagrangians on compact manifolds. We refer the reader to \cite{SorrentinoLectNotes} (and references therein) for a more detailed presentation, which include proofs of the main results.

Let $\calM(L)$ be the space of probability measures $\mu$ on 
$TM$  that are invariant  under the Euler-Lagrange flow of $L$ and such that $\int_{TM} L\,d\m<\infty$ (finite action).
We shall hereafter assume that $\calM(L)$ is endowed with the {\it vague topology}, \ie the weak$^*$--topology induced by the space $C^0_{\ell}$ of continuous functions $f:\rT M \longrightarrow \R$ having at most linear growth:
$$
\sup_{(x,v)\in \rT M} \frac{|f(x,v)|}{1+\|v\|} <+\infty\,.
$$
It is not difficult to check that $\calM(L)\subset \left(C^0_{\ell}\right)^*$. 

In the case of an autonomous Tonelli Lagrangian, it is easy to see that $\calM(L)$ is non-empty. In fact,  recall that because of the conservation of the energy
$E(x,v):=H \circ \cL(x,v) = \left\langle \frac{\dpr L }{\dpr v}(x,v),\,v \right\rangle_x - L(x,v)$ along the orbits, each energy level of $E$ is compact (it follows from  the superlinearity condition) and invariant under $\Phi^L_t$. It is a classical result in ergodic theory (sometimes called Kryloff--Bogoliouboff theorem) that a flow on a compact metric space has at least an invariant probability measure.

To each  $\m \in \calM(L)$, we may  associate  its {\it average action}:
$$ A_L(\m) = \int_{TM} L\,d\m\,. $$

The action functional 
$A_L : \calM(L) \longrightarrow \R$ is lower semicontinuous with the vague topology on $\calM(L)$ (this functional might not be necessarily continuous, see \cite[Remark 2-3.4]{contrerasiturriaga}).

In particular, this implies that there exists $\m\in\calM(L)$, which minimizes $A_L$ over $\calM(L)$.

\begin{Def}
A measure $\mu\in\calM(L)$, such that $A_L(\mu)=\min_{\calM(L)}A_L$, is called an {\it action-minimizing measure} of $L$.
\end{Def}

As we have already seen in subsection \ref{sec2.3},  by modifying the Lagrangian (without changing the Euler-Lagrange flow) one can find many other ``interesting'' measures besides those found by minimizing $A_L$. 
A similar idea can be implemented for a general Tonelli Lagrangian.
Observe, in fact,  that if $\eta$ is a $1$-form on $M$, we can interpret it as a function 
on the tangent bundle (linear on each fibre)
\beqano 
\hat{\eta}: TM &\longrightarrow&  \R \\
(x,v) &\longmapsto& \langle \eta(x),\, v\rangle_x
\eeqano
and consider a new Tonelli Lagrangian $L_{\eta}:= L - \hat{\eta}$. The 
associated Hamiltonian will be given by
$H_{\eta}(x,p) = H(x,\eta(x) + p)$.

Observe that:
\begin{itemize}
\item[i)] If $\eta$ is closed, then $L$ and $L_{\eta}$ 
have the same Euler-Lagrange flow on $T M$. See \cite{Mather91}. 
\item[ii)] If $\m \in \calM(L)$ and $\eta=df$ is an exact $1$-form, then 
$\int{\widehat{df}} d\m =0$. Thus, for a fixed $L$, the minimizing measures will 
depend only on the de Rham cohomology class $c=[\eta] \in \rH^1(M;\R)$. 
\end{itemize}

Therefore, instead of studying the action minimizing properties of a single Lagrangian, one can consider a family of such ``modified'' Lagrangians, parameterized over $H^1(M;\R)$. 
Hereafter, for any given $c\in H^1(M;\R)$, we shall denote by $\eta_c$ a closed $1$-form with that cohomology class.\\

\begin{Def}
Let $\eta_c$ be a closed $1$-form of cohomology class $c$. Then, if $\m \in \calM(L)$ minimizes 
$A_{L_{\eta_c}}$ over $\calM(L)$, we shall say that $\m$ is a {\it $c$-action minimizing measure} (or $c$-{\it minimal measure}, or {\it Mather's measure} with cohomology $c$).
\end{Def}

\noindent Compare with Fact 1 in subsection \ref{sec2.3}.\\

\begin{Rem}
Observe that the {\it cohomology class} of an action-minimizing invariant probability measure  is not intrinsic in the measure itself nor in the dynamics, but it depends on the specific choice of the  Lagrangian $L$. Changing the Lagrangian $L \longmapsto L-\eta$ by a closed $1$-form $\eta$, we shall change all the cohomology classes of its action minimizing measures by $-[\eta] \in \rH^1(M;\R)$. Compare also with Remark \ref{Remrotationvector} ({\it ii}).
\end{Rem}

One can consider the following function on $\rH^1(M;\R)$ (the {\it minus} sign is introduced for a convention that might probably become clearer later on):
\beqa{defalfa}
\a: H^1(M;\R) &\longrightarrow& \R \nn\\
c &\longmapsto& - \min_{\m\in\calM(L)} A_{L_{\eta_c}}(\m)\,.
\eeqa
This function $\alpha$ is well-defined (it does not depend on the choice of the representatives of the cohomology classes) and it is easy to see that it is convex. This is generally known as {\it Mather's $\alpha$-function}.
We have seen in subsection \ref{sec2.3} that for an integrable Hamiltonian $H(x,p)=\h( p)$, $\a(c)=\h(c)$. For this and several other reasons that we shall see later on, this function is sometimes called {\it effective Hamiltonian}. In particular, it can be proven that  $\a(c)$ is related to the energy level containing such $c$-action minimizing measures  \cite{Carneiro}.\\

We shall denote by $\calM_c(L)$ the subset of $c$-action minimizing measures:
$$\calM_c := \calM_c(L)= \{\m \in \calM(L): \; A_{L}(\m)<+\infty \;{\rm and}\; A_{L_{\eta_c}}(\m)
=-\a (c)\}.$$

We can now define a first important family  of invariant sets: 
the {\it Mather sets}. 

\begin{Def}
For a cohomology class $c \in H^1(M;\R)$, we define 
the {\it Mather set of cohomology class} $c$  as:
\be \widetilde{\cM}_c := {\bigcup_{\m \in \calM_c} {\rm supp}\,\m} 
\subset TM\,.\label{2.3}\ee
The projection on the base manifold $\cM_c = \pi \left(\widetilde{\cM}_c\right)
\subseteq M$ is called {\it projected Mather set} (with cohomology class $c$).
\end{Def}

Properties of this set:
\begin{itemize}
\item[i)] It is non-empty, compact and invariant \cite{Mather91}.
\item[ii)] It is contained in the energy level corresponding to $\alpha( c)$ \cite{Carneiro}.
\item [iii)] In \cite{Mather91} Mather proved the celebrated {\it graph theorem}: {\it let $\pi: \rT M \longrightarrow M$ denote the canonical projection. Then,
 $\pi|{\widetilde{\cM}_c}$ is an injective 
mapping of $\widetilde{\cM}_c$ into $M$, and its inverse $\pi^{-1}: \cM_c 
\longrightarrow \widetilde{\cM}_c$ is 
Lipschitz.}\\
 \end{itemize}

\vspace{10 pt}
Now, we would like to shift our attention to a related problem. As we have seen in section \ref{sec2.3}, instead of considering different minimizing problems over $\calM(L)$, obtained by modifying the Lagrangian $L$, one can alternatively try to minimize the Lagrangian $L$ putting some ``{\it constraints}'', such as, for instance, fixing the {\it rotation vector} of the measures. In order to generalize this to Tonelli Lagrangians on compact manifolds, we first need to define what we mean by rotation vector of an invariant measure.

Let $\mu\in \calM(L)$. Thanks to the superlinearity of $L$, 
the integral  $ \int_{TM} \hat{\eta} d\m$ 
is well defined and finite for any 
closed 1-form $\eta$ on $M$.
Moreover,  if $\eta$ is exact, then this integral is zero, \ie 
$ \int_{TM} \hat{\eta} d\m=0$.
Therefore, one can define a linear functional: 
\beqano
H^1(M;\R) &\longrightarrow& \R \\
c &\longmapsto& \int_{TM} \hat{\eta} d\m\,,
\eeqano
where $\eta$ is any closed $1$-form on $M$ with cohomology class $c$. By 
duality, there 
exists $\rho (\m) \in H_1(M;\R)$ such that
$$
\int_{TM} \hat{\eta} \,d\m = \langle c,\rho(\m) \rangle
\qquad \forall\,c\in H^1(M;\R)$$ 
(the bracket on the right--hand side denotes the canonical pairing between 
cohomology and 
homology). We call $\rho(\m)$ the {\it rotation vector} of $\m$. This rotation vector is 
the same as the Schwartzman's asymptotic cycle of $\mu$ (see \cite{Schwartzman} and \cite{FGS} for more details).

\begin{Rem}\label{Remrotationvector}
({\it i}) It is possible to provide a more ``geometric'' interpretation of this. Suppose for the moment that $\mu$ is ergodic. Then, it is known that a generic orbit $\g(t):=\pi \Phi^L_t(x,v)$, where $\pi:\rT M \longrightarrow M$ denotes the canonical projection, will return infinitelyf often close (as close as we like) to its initial point $\g(0)=x$. We can therefore consider a sequence of times $T_n \to +\infty$ such that $d(\g(T_n),x)\to 0$ as $n\to +\infty$, and consider the closed loops $\s_n$ obtained by ``closing'' $\g|[0,T_n]$ with the shortest geodesic connecting $\g(T_n)$ to $x$. Denoting by $[\s_n]$ the homology class of this loop, one can verify (see \cite{Schwartzman})  that $\lim_{n\to\infty} \frac{[\s_n]}{T_n} = \rho(\mu)$, independently of the chosen sequence $\{T_n\}_n$. In other words, in the case of ergodic measures, the rotation vector tells us how on average a generic orbit winds around $\rT M$. If $\mu$ is not ergodic, $\rho(\mu)$ loses this neat geometric meaning, yet it may be interpreted as the average of  the rotation vectors of its different ergodic components. 

({\it ii}) It is clear from the discussion above that the rotation vector of an invariant measure depends only on the dynamics of the system  (\ie the Euler-Lagrange flow) and not on the chosen Lagrangian. Therefore, it does not change when we modify our Lagrangian by adding a closed one form.\\
\end{Rem}

Using that the action functional $A_L: \calM(L) 
\longrightarrow \R$ is lower semicontinuous, one can prove that the map $\rho: \calM(L) \longrightarrow \rH_1(M;\R)$ is continuous and surjective, 
{\it i.e.}, for every  $h\in H_1(M;\R)$  there exists $\m\in \calM(L)$ with  $A_L(\m) < \infty$ and $\rho(\m)=h$ (see \cite{Mather91}).\\

Following Mather \cite{Mather91}, let us consider the minimal value of the average action $A_L$ over the 
probability measures with rotation vector $h$. Observe that this minimum is actually achieved because of the lower semicontinuity of $A_L$ and the compactness of $\rho^{-1}(h)$ ($\rho$ is continuous and $L$ superlinear). Let us define
\bea \label{defbeta}
\b: H_1(M;\R) &\longrightarrow& \R \nn\\
h &\longmapsto& \min_{\m\in\calM(L):\,\rho(\m)=h} A_L(\m)\,.\label{2.2}
\eea
This function $\beta$ is what is generally known as {\it Mather's 
$\beta$-function} and it is immediate to check that it is convex. 
We have seen in subsection \ref{sec2.3}, that if we have an integrable Tonelli Hamiltonian $H(x,p)=\h(p)$ and the associated Lagrangian $L(x,v)=\ell(v)$, then $\beta(h)=\ell(h)$. For this and several other reasons, this function is sometime called {\it effective Lagrangian}.\\

We can now define what we mean by action minimizing measure with a given rotation vector.

\begin{Def}
A measure $\m \in \calM(L)$ realizing the minimum in (\ref{2.2}), \ie such that $A_L(\m)~=~\b(\rho(\m))$, is called an {\it action minimizing} (or {\it minimal}, or {\it Mather's}) {\it measure} with rotation vector $\rho(\m)$.
\end{Def}

\noindent Compare with Fact 2 in subsection \ref{sec2.3}.\\

 We shall denote by $ \calM^h(L)$ the subset of action minimizing measures with rotation vector $h$:
 $$ \calM^h := \calM^h(L) = \{\m \in \calM(L): \; A_L(\m)<+\infty, \;\rho(\m)= h \; 
{\rm and}\;  A_L(\m)=\beta(h)\}.$$

 This allows us to define another important familty of invariant sets.

 \begin{Def}
 For a homology class (or rotation vector) $h\in H_1(M;\R)$, we define the
{\it Mather set corresponding to a rotation vector} $h$ as
\be \widetilde{\cM}^h := {\bigcup_{\m \in \calM^h} {\rm supp}\,\m} 
\subset TM\,,\label{2.5}\ee
and the projected one as $\cM^h = \pi \left(\widetilde{\cM}^h\right) 
\subseteq M$. 
 \end{Def}

 Similarly to what we have already seen above, this set satisfies the following properties:
\begin{itemize}
\item[i)] It is non-empty, compact and invariant.
\item[ii)] It is contained in a given energy level.
\item [iii)] It also satisfies the {\it graph theorem}: {\it let $\pi: \rT M \longrightarrow M$ denote the canonical projection. Then,
 $\pi|{\widetilde{\cM}^h}$ is an injective 
mapping of $\widetilde{\cM}^h$ into $M$, and its inverse $\pi^{-1}: \cM^h 
\longrightarrow \widetilde{\cM}^h$ is 
Lipschitz.}\\
 \end{itemize}
 
 \vspace{10 pt}

The above discussion leads to two equivalent formulations of the minimality 
of a measure $\m$:
\begin{itemize}
\item there exists a homology class $h \in H_1(M;\R)$, namely its 
rotation vector $\rho(\m)$, such that $\m$ minimizes $A_L$ amongst all 
measures in $\calM(L)$ with rotation vector $h$; \ie $A_L(\m)=\b (h)$.
\item There exists a cohomology class $c \in H^1(M;\R)$, such that $\m$ minimizes $A_{L_{\eta_c}}$ 
amongst all probability measures in $\calM(L)$; \ie $A_{L_{\eta_c}}(\m)=-\a (c)$.\\
\end{itemize}

\noindent What is the relation between these two different approaches? Are they equivalent, \ie 
$\bigcup_{h \in H_1(M;\R)} \calM^h = \bigcup_{c \in H^1(M;\R)} 
\calM_c\,$ ?\\

In order to comprehend the relation between these two families of action-minimizing measures, we need to understand better the properties of these functions $\a$ and $\beta$.
Let us start with the following trivial remark.

\begin{Rem}\label{relalphabetainteg} As we have seen in subsection \ref{sec2.3}, if we have an integrable Tonelli Hamiltonian $H(x,p)=\h(p)$ and the associated Lagrangian $L(x,v)=\ell(v)$, then $\a(c)=\h(c)$ and $\beta(h)=\ell(h)$. In this case, the cotangent bundle $\rT^*\T^n$ is foliated by invariant tori ${\mathcal T}^*_c:=\T^n\times\{c\}$ and the tangent bundle $\rT\T^n$ by invariant tori $\widetilde{\mathcal T}^h:=\T^n\times\{h\}$. In particular, we proved that 
$$
\widetilde{\cM}_c = \cL^{-1}({\mathcal T}^*_c) = \widetilde{\mathcal T}^h =  \widetilde{\cM}^h,
$$
where $\cL$ denotes the Legendre transform and $h$ and $c$ are such that  $h=\nabla \h(c) = \nabla \a(c)$ and
$c= \nabla \ell(h) = \nabla \beta(h)$.
\end{Rem}

 We would like to investigate whether a similar relation holds or does not hold in the general case. Of course, one main difficulty is that in general the {\it effective Hamiltonian} $\a$ and the {\it effective Lagrangian} $\b$, although being convex and superlinear, are not necessarily differentiable. \\
Before stating  the main relation between these two functions, let us recall some definitions and results from classical convex analysis (see \cite{Rockafellar}). Given a convex function $\varphi: V \longrightarrow \R\cup \{+\infty\}$ on a finite dimensional vector space $V$, one can consider a {\it dual} (or {\it conjugate}) function  defined on the dual space $V^*$, via the so-called {\it Fenchel transform}: $\f^*(p):= \sup_{v\in V} \big(p\cdot v - \f(v)\big)$.

\begin{Prop}\label{alphabetasuperlinear}
$\a$ and $\b$ are convex conjugate, \ie
$\a^* = \b$ and $\b^* = \a$. In particular, it follows that $\a$ and $\b$ have superlinear growth.
\end{Prop}

 Next proposition  will allow us to clearify the relation (and duality) between the two minimizing procedures above. To state it, recall that, like any convex 
function on a
finite-dimensional space, $\b$
admits a subderivative at each point $h\in \rH_1(M;\R)$, \ie we can find $c\in 
\rH^1(M;\R)$ such that
\beqa{subdiff}
\forall h'\in \rH_1(M;\R), \quad \b(h')-\b(h)\geq \langle c,h'-h\rangle.\eeqa
As it is usually done, we shall denote by $\partial \b(h)$ the set of $c\in 
\rH^1(M;\R)$ that 
are subderivatives of $\b$ at $h$, \ie the set of $c$ which satisfy the 
inequality above.  Similarly, we shall denote by $\dpr \a(c)$ the set of subderivatives of $\a$ at $c$.\\
Fenchel's duality implies an easier characterization of subdifferentials:
{ \it $c\in \partial \b(h)$ if and only if $\langle c,h\rangle=\a(c )+\b(h)$} (similarly for $h\in \partial \a©$).\\

We can now prove that what was observed in Remark \ref{relalphabetainteg} continues to hold in the general case.

\begin{Prop}\label{CaracMinim}
Let  $\m \in \calM(L)$ be an invariant probability measure. Then:\\
{\rm (i)} $A_L(\m)=\b(\rho(\m))$ {if and only if}  there exists $c\in H^1(M;\R)$ such that $\m$ minimizes $A_{L_{\eta_c}}$  
{\rm(}\ie $A_{L_{\eta_c}}(\mu)=-\a(c)${\rm)}.\\
{\rm (ii)} If $\m$ satisfies $A_L(\m)=\b(\rho(\m))$ and $c \in 
\rH^1(M;\R)$, then $\m$ minimizes $A_{L_{\eta_c}}$ if and only if $c\in \partial 
\b(\rho(\mu))$ {\rm(}or equivalently $\langle c,h\rangle=\a(c)+\b(\rho(\m)\rm{)}$.
\end{Prop}

\begin{Rem}\label{remarkinclusionsmathersets}
({\it i}) It follows from the above proposition that both minimizing procedures lead to the same sets of invariant probability measures:
$$
\bigcup_{h \in H_1(M;\R)} \calM^h = \bigcup_{c \in H^1(M;\R)} 
\calM_c\,.
$$
In other words,  minimizing over the set of invariant measures with a fixed rotation vector or minimizing - globally - the modified Lagrangian (corresponding to a certain cohomology class) are dual problems,  as the ones that often appear in linear programming and optimization.

({\it ii}) In particular, we have the following inclusions between Mather sets:
$$
 c\in \dpr \beta(h) \quad \Longleftrightarrow \quad h \in \partial \a(c) \quad \Longleftrightarrow \quad \widetilde{\cM}^h \subseteq \widetilde{\cM}_c\,.
 $$
Moreover, for any $c\in \rH^1(M;\R)$:
\begin{equation}\label{inclusionMathersets}
\widetilde{\cM}_c = \bigcup_{h\in \dpr \a(c)} \widetilde{\cM}^h\,.
\end{equation}
\end{Rem}

\vspace{10 pt}

\begin{Rem}\label{rem2.17}
({\it i}) In the above discussion we have only discussed properties of invariant probability measures associated to the system. Actually, one could study directly orbits of the systems and look for orbits that {\it globally minimize} the action of a modified Lagrangian (in the same spirit as before). This would lead to the definition of two other families of invariant compact sets, the {\it Aubry sets} $\widetilde{\cA}_c$  and the {\it Ma\~n\'e sets} $\widetilde{\cN}_c$, which is also parameterized by $\rH^1(M;\R)$ (the parameter which describes the modification of the Lagrangian, exactly in the same way  as before). For a given $c\in \rH^1(M;\R)$, these sets contain  the Mather set $\widetilde{\cM}_c$, and this inclusion may be strict. In fact, while the motion on the Mather sets is {\it recurrent} (it is the union of the supports of invariant probability measures), the Aubry and the Ma\'n\'e sets may contain non-recurrent orbits as well.\\
({\it ii }) Differently from what happens with invariant probability measures, it will not be always possible to find {\it action-minimizing orbits} for any given rotation vector (not even possible define a rotation vector for every action minimizing orbit). For instance, an example due to Hedlund \cite{Hedlund} provides the existence of a Riemannian metric on a three-dimensional torus, for which minimal geodesics exist only in three directions. This can be extended to any dimension larger than three.\\
\end{Rem}



\section{Weak Liouville--\Arnold theorem and its implications}\label{sec3}

In this section we would like to discuss the  main results stated in Section \ref{sec1.1}.
Roughly speaking, the main idea behind our approach consists in studying how the existence of independent integrals of motion of a Tonelli Hamiltonian $H$ relates to the structure and ``size'' of its Mather sets. Moreover, using the {\it symplectic invariance} of these sets, one is  be able to recover the involution hypothesis at least {\it locally}. Hereafter we shall mainly work in the Hamiltonian setting. Let us denote by $\cM^*_c$ the corresponding Mather set in the cotangent bundle, {\it i.e.},
$\cM^*_c = \cL (\widetilde{\cM}_c)$, where $\cL$ denotes the Fenchel-Legendre transform.\\

The key properties that we use, can be summarized as follows:
\begin{itemize}
\item[P1 -] {the Mather  sets are invariant under the flow of any integral of motion of $H$};
\item[P2 -] {the existence of $k$ independent integrals of motion implies that the ``size'' of each Mather set is bigger or equal than $k$} (in a sense that will be explained below); 
\item[P3 -] {the integrals of motion are locally in involution on the Mather sets}.\\
 \end{itemize}

Let us discuss these properties more in details.\\

\noindent{\bf P1 - Symplectic invariance of Mather sets.}  In \cite{SorrentinoTAMS} I proved the following result:

\begin{Prop} \label{maintheorem}
Let $H$ be a Tonelli Hamiltonian on $\rT^*M$ and $F$ an integral of motion of $H$.  
Let us denote by $\Phi_H$ and $\Phi_F$ the respective flows. Then, the following holds:
\begin{itemize}
\item[{\rm (i)}] If $\m$ is a $c$-action minimizing measure of $H$, then ${\Phi^t_F}_*\m$ is still a $c$-action minimizing measure of $H$, for each $t\in \R$,
where the lower $*$ denotes the push-forward of the measure. 
\item[{\rm (ii)}] The Mather set $\cM^*_c$ is invariant under the action of $\Phi^t_F$, for each $t\in \R$ and for each $c\in \rH^1(M;\R)$. In particular, for each $t\in \R$, $\Phi^t_F$  maps each connected component of  $\cM^*_c$ into itself.
 \end{itemize}
\end{Prop}

\vspace{10 pt}

\begin{Rem}\label{giadim}
It is worthwhile to point out that this result can be also deduced from a  
result by Patrick Bernard  \cite[Theorem in Section $1.10$, page $6$]{Bernard} on the symplectic invariance of the Mather and Aubry sets.  
In fact for any fixed time $t$ the Hamiltonian flow $\Phi_F^t$ is an exact symplectomorphism that preserves $H$. \\
Another related result is contained in \cite{Maderna}, where the author considers the action of symmetries of the Hamiltonian, \ie $C^1$-diffeomorphisms of $M$ that preserve $H$. One can deduce from the results therein that the Mather and Aubry sets of $H$ are invariant under the action of the  connected component of the identity in the group of such diffeomorphisms. From our point of view, these diffeomorphisms correspond to integrals of motion depending only on the $x$-variables.\\
\end{Rem}

\noindent{\bf P2 - Independent Integrals of motion and size of the Mather sets.}
As recalled in Section \ref{sec1.1}, Liouville--\Arnold theorem is concerned with {\it independent} integrals of motion, \ie integrals of motion whose differentials are linearly independent, as vectors, at each point.  Let us see how the existence of independent integrals of motion relates to  the ``size'' of the Mather  of $H$. 
In order to make clear what we mean by ``size'' of these sets, let us introduce some notion of tangent space.
We shall call {\it generalized tangent space} to $\cM^*_c$ at a point $(x,p)$, the set of all vectors that are tangent to curves in $\cM^*_c$  at $(x,p)$. We shall denote it by $\rT_{(x,p)}\cM^*_c$  and we shall define its {\it rank} to be the largest number of linearly independent vectors that it contains.
 In particular, if the Mather set does not contain any fixed point
(\ie $dH(x,p)\neq 0$ for all $(x,p)\in \cM^*_c$), then $ {\rm rank}\ \rT_{(x,p)}\cM^*_c \geq 1$; 
in fact, since these sets are invariant, the Hamiltonian vector field $X_H(x,p) \neq 0$ is tangent to them. 

\begin{Prop}\label{Prop100}
Let $H$ be a Tonelli Hamiltonian on $\rT^*M$ and suppose that there exist $k$ independent integrals of motion.
Then, ${\rm rank}\ \rT_{(x,p)}\cM^*_c \geq k$ at all points $(x,p)\in \cM^*_c$ for each $c\in \rH^1(M;\R)$.\\
\end{Prop}

The proof of this result follows easily from Proposition \ref{maintheorem}. More specifically, it follows from the fact that $\cM^*_c$ is invariant under the flows of the $k$ independent integrals of motion. The linear independence of the corresponding vector fields (which are therefore tangent to this set) follows from the independence of the integrals of motion  and the non-degeneracy of the symplectic form $\omega$.\\

 In particular, the existence of the maximum possible number of integrals of motion (\ie $k=n$) implies that these sets are invariant smooth Lagrangian graphs (see \cite[Remark 3.5 and Lemmata 3.4 \& 3.6]{SorrentinoTAMS}).\\ 
 
\begin{Cor}\label{weakLiouville}
Let $H$ be a weakly integrable Tonelli Hamiltonian on $\rT^*M$.  Then, for each $c\in \rH^1(M;\R)$ such that 
$\cM^*_c \subset {\rm reg}(F)$, then  we have that $\cM^*_c$ projects over the whole $M$ and therefore it is an invariant Lipschitz Lagrangian graph.
\end{Cor}

 In particular, smoothness is a consequence of the fact that these graphs lie in level sets of the integral map, which is non-degenerate (see 
\cite[Lemma 3.6]{SorrentinoTAMS}). \\

\noindent{\bf P3 - Local involution of the Mather sets.}

The most important peculiarity of these action-minimizing sets, firstly observed in \cite{SorrentinoTAMS}, is that they force the integrals of motion to Poisson-commute on them. In fact,  one can recover the involution property of the integrals of motion, at least locally (see \cite[Proposition 27]{SorrentinoTAMS} for a detailed proof).

\begin{Prop}\label{prop2Sor09}
Let $H$ be a Tonelli Hamiltonian on $T^*M$ and let $F_1$ and $F_2$ be two integrals of motion. Then for each $c\in H^1(M;\R)$ we have that $\{F_1,F_2\}(x,\hat{\pi}_c^{-1}\!(x))=0$ for all $x\in \overline{{\rm Int}\big(\cM_c\big)}$, where $\hat{\pi}_c=\pi|\cM^*_c$   and $\cM_c=\pi\big(\cM^*_c\big)$.
\end{Prop}

\begin{Rem}
Observe that the above set $\overline{{\rm Int}\big(\cM_c\big)}$ may be empty. What the proposition says is that whenever it is non-empty, the integrals of motion are forced to Poisson-commute on it. In the cases that we shall be considering hereafter, $\cM_c=M$ -- since it corresponds to the projection of a Lagrangian graph -- and therefore it is not empty.\\
\end{Rem}

\noindent{\bf Sketch of the proof of Theorem \ref{mainthm} ( Weak Liouville-\Arnold Theorem)}. 
{\it i}) The existence of these smooth Lagrangian graphs follows from Corollary \ref{weakLiouville}. 
In fact, since $\cM^*_c$ is contained in the set of regular points of $F$,  it follows that the Mather set $\cM_c^*$ is a $C^1$ invariant Lagrangian graph $\Lambda_c$ of cohomology class $c$. Therefore, $\Lambda_c$ supports an invariant probability measure of full support. Using upper semi-continuity\footnote{Upper semicontinuity of the Mather sets does not hold in general. However, in this case it does, since it coincides with a Lipschitz Lagrangian graph (see for instance \cite[Proposition 13]{Arn08}).} of the Mather sets 
to deduce that the Mather sets corresponding to nearby cohomology classes must also lie in ${\rm reg}\ {F}$. Hence, there  exists an open neighbourhood $\cO$ of $c$ in $H^1(M;\R)$ such that  $\cM_{c'}^*  \subset {\rm reg}\ {F}$ for all $c'\in \cO$ and applying the same arguments as above, we can conclude that each $\cM_{c'}^*$ is a smooth invariant Lagrangian graph of cohomology class $c'$.\\
For the proof of the fact that such Lagrangian graphs admits the structure of a smooth $\T^d$-bundle over a parallelisable base $B^{n-d}$, for some $d>0$, we refer the reader to \cite[Proposition 2.7]{ButlerSorrentino}.\\
{\it ii}) It follows from above that these $\Lambda_{c'} = \cM^*_{c'}$, therefore they coincide with the union of supports of invariant probability measures. The fact that these measures have all the same rotation vector and that the orbits in their supports are conjugate by a smooth diffeomorphism isotopic to the identity, is discussed in \cite[Remark 3.1]{ButlerSorrentino}.\\
{\it iii}) We know from ({\it ii}) that these graphs are Schwartzman uniquely ergodic, \ie all invariant probability measures on  $\Lambda_{c'}$ have  the same rotation vector $h_{c'} \in H_1(M;\R)$. 
The differentiability of $\alpha$ follows then from  \cite[Corollary 3.6]{FGS}. The differentiability of $\beta$ follows the disjointness of these graphs (see for instance \cite[Theorem 3.3]{FGS} or \cite[Remark 4.26 (ii)]{SorrentinoLectNotes}).\\

Now we discuss whether or not there are cases in which this weaker notion of integrability is equivalent to the classical  one (in the sense of Liouville). As remarked in Section \ref{sec1.1}, the union of  these Lagrangian graphs is not necessarily a foliation of the whole phase space. In fact, if the dimension of $H^1(M;\R)$ is less than the dimension of $M$, this family of graphs is not sufficient to foliate $\rT^*M$ or even to have non-empty interior (for instance, think about the case in which $\rH^1(M;\R)$ is trivial). What we prove in Theorem    \ref{maincor} is that when this obstacle is removed, then the two notions coincide.\\

\noindent{\bf Sketch of the proof of Theorem \ref{maincor}}. 
Let us denote $\Lambda_{c'}=\{(x,\lambda_{c'}(x)):\; x\in M\}$ as usual. Observe that the map:
\begin{eqnarray*}
\Psi: \cO \times M &\longrightarrow & T^*M \\
(c', x) &\longmapsto& \lambda_{c'}(x)
\end{eqnarray*}
is continuous (see \cite[proof of Theorem 1.2 (iii)]{ButlerSorrentino}). 
If $\dim H^1(M;\R) \geq \dim M$, then the continuity of $\Psi$ implies that these Lagrangian graphs 
$\Lambda_{c'}$ foliate an open neighbourhood of $\Lambda_c$. It follows from 
Proposition \ref{prop2Sor09} that the components of $F$ commute in this open region. Therefore, each $\Lambda_{c'}$ is an $n$-dimensional manifold which is invariant under the action of $n$ commuting vector fields, which are linearly independent at each point. It is a classical result that $\Lambda_{c'}$ must be then diffeomorphic to an $n$-dimensional torus and that the motion on it is conjugate to a rotation (see for instance \cite{Arnoldbook}). \\


\section{Minimal average action and Integrability}\label{sec4}
In Section \ref{sec2.4} we have introduced the {\it minimal average action(s)}, namely the so-called Mather's $\alpha$ and $\beta$ functions:
$$\a: \rH^1(M;\R) \longrightarrow \R \quad {\rm and}\quad  \b:\rH_1(M;\R) \longrightarrow \R,$$
defined respectively in (\ref{defalfa}) and (\ref{defbeta}).

In Section \ref{sec2} (more specifically in subsections \ref{sec2.3} and \ref{sec2.4}) we have discussed the following results, which relate the differentiability of these functions (or lack thereof) to properties of the Mather sets:
\begin{itemize}
\item[{\it i)}]  if we have an integrable Tonelli Hamiltonian $H(x,p)=\h(p)$ and the associated Lagrangian $L(x,v)=\ell(v)$, then $\a(c)=\h(c)$ and $\beta(h)=\ell(h)$. In particular, we proved that  $$
\widetilde{\cM}_c =   \widetilde{\cM}^h,
$$
where $h$ and $c$ are such that  $h= \nabla \a(c)$ and $c=  \nabla \beta(h)$. See Remark \ref{relalphabetainteg}. 
\item[{\it ii)}]
We have the following inclusions between Mather sets:
$$
 c\in \dpr \beta(h) \quad \Longleftrightarrow \quad h \in \partial \a(c) \quad \Longleftrightarrow \quad \widetilde{\cM}^h \subseteq \widetilde{\cM}_c\,.
 $$
\item[{\it iii)}] For any $c\in \rH^1(M;\R)$:
$$
\widetilde{\cM}_c = \bigcup_{h\in \dpr \a(c)} \widetilde{\cM}^h\,.
$$
See Remark \ref{remarkinclusionsmathersets}.\\
\end{itemize}

In particular, we can deduce immediately that:
\begin{itemize}
\item lack of differentiability of $\alpha$ at $c$  $\Longrightarrow$  the Mather set $\widetilde{\cM}_c$ contains action-minimizing measures with different rotation vectors.
\item lack of differentiability of $\beta$ at $h$  $\Longrightarrow$  two Mather sets $\widetilde{\cM}_c$ and $\widetilde{\cM}_{c'}$  have non-trivial intersection (they must both contain $\widetilde{\cM}^h$).\\
\end{itemize}

Hence, a necessary condition for the system to be completely integrable, is that $\alpha$ and $\beta$ are $C^1$.
One can weaken the assumption on the complete integrability of the system and consider  $C^0$-{\it integrable systems} (this notion has been first introduced by Marie-Claude Arnaud in \cite{Arnaud}):\\

\begin{Def}\label{C0integrability}
A Tonelli Hamiltonian $H:\rT^*M \longrightarrow \R$ is said to be $C^0$-integrable, if there exists a foliation of $\rT^* M$ made by invariant Lipschitz Lagrangian graphs, one for each cohomology class.\\
\end{Def}

\begin{Rem}\label{remmmm}
({\it i})
Let $M$ be a compact manifold of any dimension, $L: \rT M \longrightarrow \R$ a Tonelli Lagrangian and $H: \rT^* M \longrightarrow \R$ the associated Hamiltonian. If $H$ is $C^0$-integrable, then $\beta$ is $C^1$.  See \cite[Lemma 5]{MassartSorrentino}.\\
({\it ii}) As we proved in Theorem \ref{mainthm} (iii), the notion of weak integrability also implies the  differentiability of $\beta$ and $\alpha$ in an open set.\\
({\it iii}) The notion of $C^0$-integrability is conceptually weaker than the notions of Liouville integrability and Weak integrability. It is an open problem whether there exist systems that are $C^0$ integrable, but not Liouville integrable or weakly integrable.
\end{Rem}

These observations raise the following question, which is the starting point of the work in \cite{MassartSorrentino}: {\it is the converse true?} Namely, {\it  is true that if $\beta$ is $C^1$ then the system is $C^0$-integrable?}\\

The question stated in this form has clearly a negative answer. Examples of Lagrangians  admitting a smooth $\beta$-function, but not integrable, are easy to construct: 
\begin{itemize}
\item[-] trivially, if the base manifold $M$ is such that $\dim \rH_1(M;\R)=0$ then $\beta$ is a function
defined on a single-point set and it is therefore smooth. 
\item[-] if $\dim \rH_1(M;\R)=1$ then a result by Carneiro \cite{Carneiro} (namely, $\beta$ is always differentiable in the radial direction, see Lemma \ref{Carneiro radial}) allows one to conclude that $\beta$ is differentiable everywhere, except possibly at the origin. \\
\end{itemize}

Therefore, one should rephrase the question in the following way:\\

\noindent{\bf Question:}  {\it With the exception of the mentioned trivial cases (i.e, when $\dim \rH_1(M;\R)\leq 1$),
does the regularity of $\b$ imply the integrability of the system?}\\

In \cite{MassartSorrentino} we addressed this question in the case of closed surfaces, not necessarily orientable (in this latter case, one considers the lifted Lagrangian to the orientable double cover).  \\

Let us start by recalling some terminology:  
\begin{itemize}
\item a homology class $h$ is said to be $k$-irrational, if $k$ is the dimension of the smallest subspace of $\rH_1(M;\R)$ generated by integer classes and containing $h$. In particular, $1$-{\it irrational} means ``on a line with rational slope'', while {\it completely irrational} stands for ``$\dim \rH_1(M;\R)$-irrational''. 
\item  A homology $h$ is said to be {\it singular} if its Legendre transform $\dpr \beta(h)$ is a {\it singular flat}, \ie its Mather set $\widetilde{\cM}(\dpr \beta(h))$   contains fixed points. Observe that the set of singular classes, unless it is empty,  contains the zero class and  is compact.
\item For $h \in \rH_1(M;\R)\setminus \left\{0\right\}$, we define the {\it maximal radial flat} $R_h$ of $\beta$ containing $h$ as the largest subset of the half-line $\left\{th \; : t \in \left[0,+\infty \right) \right\}$ containing $h$ (not necessarily in its relative interior) in restriction to which $\beta$ is affine. 
\end{itemize}

In \cite{Massart} Daniel Massart proved the following result.

\begin{Teo}[{\bf Massart \cite[Theorem 3]{Massart}}]
Let $M$ be a closed surface and $L$ be an autonomous Tonelli Lagrangian on $\rT M$. If 
$h$ is a $1$-irrational, nonsingular homology class, then $\widetilde{\cM}^h$  is a union of periodic orbits.
\end{Teo}

\begin{Rem}
This results does not hold anymore in dimension greater than two.\\
\end{Rem}

\subsection{Differentiability of Mather's $\beta$--function on closed surfaces}
Let us now discuss some differentiability properties of Mather's $\b$--function at $1$-irrational homology classes. These will play a crucial role in the proof of the main results stated in subsection \ref{sec1.2}.
 First let us get the Klein bottle case out of the way. We shall use the following lemma from \cite{Carneiro}:

\begin{Lem}\label{Carneiro radial}
If $L$ is an autonomous Tonelli Lagrangian on a closed manifold $M$, then at every $h \in \rH_1(M;\R) \setminus \{0\}$, $\beta$ is differentiable in the radial direction, that is, the map 
$$
\begin{array}{rcl}
B_h: \R & \longrightarrow & \R \\
t & \longmapsto & \beta (th)
\end{array}
$$
is $C^1$ on $\R\setminus\{0\}$.\\
\end{Lem}

\begin{Cor}\label{Klein}
If $L$ is an autonomous Tonelli Lagrangian on the Klein bottle, then $\beta$ is $C^1$, except possibly at $0$.
\end{Cor}

\begin{Proof}
The first Betti number of the Klein bottle is one, that is, there exists $h_0 \in \rH_1(\K ;\R) \setminus \{0\}$ such that for all $h \in \rH_1(\K ;\R) $, there exists $t \in \R$ such that $h=th_0$. Then, we use Lemma \ref{Carneiro radial}.
\end{Proof}

The meaning of the next theorem is that for an autonomous Lagrangian on a closed surface, at a 1-irrational, non-singular homology class, $\beta$ is differentiable only in the directions where it is flat, and in the radial direction. Indeed, in the statement below, $\mathcal{V}(h_0)$ may be viewed as a measure of the non-differentiability of $\beta$ at $h_0$,  while $\partial \alpha (c_0)$ is the largest flat containing $h_0$ in its relative interior. See \cite[Theorem 2]{MassartSorrentino} for a complete proof.

\begin{Teo}\label{betadiff, thm}
Let  $M$ be a closed surface and $L$ an autonomous  Tonelli Lagrangian on $\rT M$. If:
\begin{itemize}
	\item $h_0$ is a 1-irrational, non-singular homology class,
	\item $(\gamma_i, \dot\gamma_i)_{i \in I}$ are the periodic orbits which comprise the supports of all action-minimizing measures with rotation vectors $th$, for all $th$ in $R_{h_0}$,
	\item $c_0$ is a cohomology class in the relative interior of $\partial \beta (h_0)$,
	\item $\mathcal{V}(h_0)$ is the vector subspace of $\rH^1(M;\R)$ generated by the differences $c_1 - c_2$, where $c_1,c_2$ are elements of 
	$\partial \beta (h_0)$,
\end{itemize}
then we have:
\begin{itemize}
  \item either 
  $$
  \mathcal{V}(h_0) = \partial \alpha (c_0)^{\perp}= \bigcap_{i \in J } h_i ^{\perp}
  $$
  where orthogonality is meant with respect to the duality between $\rH_1(M;\R)$ and $\rH^1(M;\R)$
  \item or $M= \T^2$   and  the closed curves $\gamma_i $, $i \in I$, foliate $M$; in this case $ \mathcal{V}(h_0) = \{0\}$. 
\end{itemize}
\end{Teo}

From Theorem \ref{betadiff, thm} we deduce that when $M$ is oriented, $\beta$ is never differentiable at any 1-irrational, non-singular homology class, unless $M= \T^2$ and $M$ is foliated by periodic orbits. This proves Theorem \ref{teoms} (i) (orientable case).\\

 \begin{Rem}
When $M$ is not orientable, the situation is different because $\beta$ may have flats of maximal dimension (that is, of dimension equal to the first Betti number of $M$). So $\beta$ may well be differentiable at some 1-irrational, non-singular homology class. 
For instance, by \cite[Theorem 1.3]{nonor},  there exists a Riemannian metric on $M$, whose stable norm has a (finite) polyhedron as its unit ball. Let $L$ be the Lagrangian induced by this Riemannian metric, then the $\beta$-function of $L$ is half the square of the stable norm. In particular $\beta$ is differentiable everywhere but along a finite number of straight lines. So $\beta$ is differentiable at most 1-irrational classes. Furthermore, since the Lagrangian is a Riemannian metric:
\begin{itemize}
  \item no homology class other than zero is singular
  \item the $\beta$-function is quadratic (\ie 2-homogeneous) so radial faces are trivial, \ie $\forall h, \  R_h = \{h\}$. 
\end{itemize} 
On the other hand, for every homology class the Mather set is a finite union of closed geodesics, so its projection can never be the whole manifold $M$. \\
  \end{Rem}

In the non-orientable case, in general, all we have is the following.

  \begin{Prop}\label{maintheo2, nonor}
Let $M$ be a closed, non-orientable surface other than the projective plane or the Klein bottle,  and let $L:\rT M \longrightarrow \R$ be a Tonelli Lagrangian. Then, there exists some homology class $h$ such that $\beta$ is not differentiable at $h$.
\end{Prop}

\noindent See \cite[Corollary 2]{MassartSorrentino}.  This proves Theorem \ref{teoms} (i) (non-orientable case).\\

\subsection{$C^0$--Integrability on closed surfaces} \label{sec4.2}

As defined in Definition \ref{C0integrability}, recall that  a Tonelli Hamiltonian $H:\rT^*M \longrightarrow \R$ is said to be $C^0$-integrable, if there exists a foliation of $\rT^* M$ made by invariant Lipschitz Lagrangian graphs, one for each cohomology class.

In \cite[Proposition 4]{MassartSorrentino},  we proved the following result.

\begin{Prop} \label{torus_only_integrable}
The torus is the only closed surface which admits a $C^0$-integrable Hamiltonian.\\
\end{Prop}

\noindent Observe that this proves Theorem \ref{teoms} (ii).\\

\noindent {\bf Ideas of the proof.} 
\begin{itemize}
\item First,  no Hamiltonian on the sphere can be $C^0$-integrable. Indeed, any Lagrangian graph is exact since the sphere is simply connected, and any two exact Lagrangian graphs intersect, because any $C^1$ function on the sphere has critical points. 
\item Likewise, no Hamiltonian on the projective plane can be 
$C^0$-integrable, otherwise its lift to the sphere would be $C^0$-integrable.
\item  Let $\K$ denote the Klein bottle.
 For each $x \in \K$, let us define
\beqano
F_{x} : \rH^1(\K;\R)\simeq \R &\longrightarrow & \rT^*_x \K \simeq \R^2\\
c &\longmapsto& \L_c \cap \rT^*_x \K,
\eeqano
where $\L_c$, for  $c\in\rH^1(\K;\R)$, are the Lagrangian graphs foliating $\rT^*\K$.
It is possible to check that $F_{x}$ is continuous (see \cite[Lemme 4.22]{Arnaud}) and  injective (as it follows from the disjointness of the $\L_c$'s). Moreover, if the Hamiltonian is $C^0$-integrable, the map $F_x$ is surjective. Now there is no such thing as a continuous bijection from $\R$ to $\R^2$, so there is no $C^0$-integrable Hamiltonian on the Klein bottle. 
\item  The same argument can be used for any surface with first Betti number $>2$. 
\item Finally,
no Hamiltonian on the connected sum of three projective planes can be $C^0$-integrable, otherwise it would  lift to  a $C^0$-integrable Hamiltonian on a surface of genus two.
\end{itemize}
\qed

Finally, we can prove the part (iii) of Theorem \ref{teoms}.

\begin{Teo}\label{teo3massor}
Let  $L:\rT \T^2 \longrightarrow \R$ be a Tonelli Lagrangian on the two-torus. 
Then, $\beta$ is $C^1$ if and only if the system is $C^0$-integrable.\\
\end{Teo}

\noindent {\bf Sketch of the proof.}\footnote{See \cite[Theorem 3]{MassartSorrentino} for a complete proof.}
[$\Longleftarrow$] See Remark \ref{remmmm} ({\it i}). 
 [$\Longrightarrow$] For each homology class $h$, let us denote by $c_h:=\dpr \beta(h)$. 
 If $h$ is non-singular and $1$-irrational, then  Theorem \ref{betadiff, thm} says that $\L_{c_h}:= \cM^*_{c_h}$ is an invariant Lipschitz Lagrangian graph of cohomology class $c_h$, which is foliated by periodic orbits of homology $h$ and same minimal period $T_h$. Since the dynamics on these Lipschitz Lagrangian graph is totally periodic, \ie $\Phi^H_{T_h}\big| \L_{c_h} = {\rm Id}\big| \L_{c_h}$, then it follows from the result in \cite{Arnaud} (see for instance the proof of Th\'eor\`eme 4) that this graph is in fact $C^1$.
 
One can show  that such cohomology classes $c_h$ are dense in $\rH^1(\T^2;\R)$.  Using the semicontinuity of these action-minimizing sets, we can deduce that for each $c\in \rH^1(\T^2;\R)$ there exists an invariant Lipschitz Lagrangian graph with cohomology class $c$, which we denote $\L_c$. Observe that all these $\L_c$'s are disjoint (it is a straightforward consequence of the differentiability of $\beta$).  For each ${x_0}\in \T^2$, let us define
\beqa{function}
F_{x_0} : \rH^1(\T^2;\R)\simeq \R^2 &\longrightarrow & \rT^*_{x_0}\T^2 \simeq \R^2\\
c &\longmapsto& \Lambda_c \cap \rT^*_{x_0}\T^2.\nonumber
\eeqa
This map is injective (as it follows from the disjointness of the $\L_c$'s). Moreover, one can also prove that  $F_{x_0}$ is also continuous. 
Therefore, we can conclude that $F_{x_0}(\R^2)$ is open (see for instance \cite{Brouwer}). In a similar way one can show that this image is also closed and hence that it is all of $\R^2$. Since this holds for all $x_0\in \T^2$, we can conclude that $\bigcup_c \L_c = \rT^*\T^2$, that is, the system is $C^0$-integrable. \qed\\

Summarizing, the  proof of Theorem \ref{teoms} is:
\begin{itemize}
\item[-]   the orientable case of ({\it i}) follows from Theorem \ref{betadiff, thm} and the non-orientable one from Proposition \ref{maintheo2, nonor};
\item[-] ({\it ii})  follows from Proposition \ref{torus_only_integrable};
\item[-] ({\it iii}) follows from Theorem \ref{teo3massor}.\\
\end{itemize}

Observe that in \cite{MassartSorrentino} we were able to deduce more information about the dynamics of a $C^0$-integrable system (see also \cite[Corollary 3]{MassartSorrentino} and \cite{Arnaud}).

\begin{Cor}\label{corrr}
Let $H: \rT^* M \longrightarrow \R$ be a $C^0$-integrable Hamiltonian on a two-dimensional closed manifold $M$. Then, $M$ is diffeomorphic to $\T^2$. Moreover :
\begin{itemize}
\item[{\rm (}i{\rm)}]  for each $1$-irrational homology class $h$, there exists an invariant Lagrangian graph foliated by periodic orbits with homology $h$ and the same minimal period;
\item[{\rm (}ii{\rm)}] for each completely irrational homology class, there exists an invariant Lagrangian graph on which the motion is conjugate to an irrational rotation on the torus or to a Denjoy type homeomorphism; 
\item[{\rm (}iii{\rm)}] there exists a dense $G_{\delta}$ set of (co)homology classes, for which the motion on the corresponding invariant torus is conjugate to a rotation;
\item[{\rm (}iv{\rm)}] as for the $0$-homology class, there exists a $C^1$ invariant torus  $\L_{c(0)}= \{(x,\frac{\dpr L}{\dpr v}(x,0):\; x\in \T^2\}$ consisting of fixed points.
\end{itemize}
\end{Cor}

\vspace{10 pt}

The above results proves the $C^0$--integrable systems. An open question is if the differentiability of Mather's $\beta$ function implies stronger properties, for example the integrability in the sense of Liouville.
In the case of mechanical systems we can bridge this gap, using Burago and Ivanov's theorem on metrics without conjugate points \cite{BI}.
Recall that a mechanical  Lagrangian is $L(x,v) = 1/2 \,  g_x(v,v) +f(x)$, where $g$ is a Riemannian metric and $f$ is a $C^2$ function on $\T^2$ (see also subsection \ref{sec2.2}).\\

\noindent {\bf Proposition \ref{propmech}.}
{\it Let $L$ be a  mechanical Lagrangian on a $2$-dimensional torus, whose $\beta$-function is $C^1$. Then the potential $f$ is identically constant and the metric $g$ is flat. In particular, $L$ is integrable in the sense of Liouville.}\\

\noindent {\bf Sketch of the proof.}\footnote{See \cite[Proposition 6]{MassartSorrentino} for a complete proof.} 
In the mechanical case, the only fixed points of the Euler-Lagrange flow are the critical points of the potential $f$, and the only minimizing fixed points are the minima of $f$. Hence, if $f$ is not constant, then the Lagrangian cannot be $C^0$-integrable (see Corollary \ref{corrr} (iv)). Furthermore, since the Lagrangian is $C^0$-integrable, every orbit is minimizing, in particular, it can be proven that there are no conjugate points. So by Burago and Ivanov's proof of the Hopf Conjecture \cite{BI}, the metric $g$ is flat. This completes the proof.



\end{document}